\documentclass[aos,MSNbibl,rotating,nameyear,dvips]{arximspdf}
\usepackage{mathrsfs}
\usepackage{dcolumn}
\usepackage{graphicx}

\doi{10.1214/14-AOS1232} 
\volume{42}
\issue{4}
\pubyear{2014}
\firstpage{1426}
\lastpage{1451}

\makeatletter
\newproclaim{example}{Example}
\newcolumntype{d}[1]{D{.}{.}{#1}}

\newcommand{\rrVert}{\Vert}
\newcommand{\llVert}{\Vert}
\newcommand{\eqref}[1]{(\ref{#1})}
\def\Mb{\mathbf{M}}
\def\Ib{\mathbf{I}}
\def\Fb{\mathbf{F}}
\def\cb{\mathbf{c}}
\def\ub{\mathbf{u}}
\def\wb{\mathbf{w}}
\def\xb{\mathbf{x}}
\def\yb{\mathbf{y}}
\def\0b{\mathbf{0}}
\def\SM{\mathcal{M}}

\def\D{\mathrm{d}}
\def\mt{\theta}
\def\mtb{\bar\theta}
\def\mve{\varepsilon}
\def\Ex{\mathbb{E}}
\def\mth{\hat{\theta}}
\def\ma{\alpha}
\def\me{\epsilon}
\def\ra{\rightarrow}
\def\SG{{\mathcal G}}
\def\ms{\sigma}
\def\ml{\lambda}
\def\mp{\partial}
\def\SE{\mathscr{E}}
\def\SB{\mathscr{B}}
\def\SG{\mathscr{G}}
\def\SL{\mathscr{L}}
\def\SM{\mathscr{M}}
\def\SX{\mathscr{X}}
\def\Var{\operatorname{Var}}
\def\cov{\operatorname{cov}}
\def\inter{\operatorname{int}}
\def\rank{\operatorname{rank}}
\def\dim{\operatorname{dim}}
\def\SSS{\mathbb{S}}

\newtheorem{theo}{Theorem}

\makeatother

\begin{document}
\begin{frontmatter}

\title{Optimum design accounting for the global nonlinear behavior of
the model}
\runtitle{Optimum design accounting for model nonlinearity}

\begin{aug}
\author[a]{\fnms{Andrej} \snm{P\'azman}\ead[label=e1]{pazman@fmph.uniba.sk}\thanksref{t1}}
\and
\author[b]{\fnms{Luc} \snm{Pronzato}\corref{}\ead[label=e2]{pronzato@i3s.unice.fr}}
\runauthor{A. P\'azman and L. Pronzato}
\thankstext{t1}{Supported by the VEGA Grant No. 1/0163/13.}
\affiliation{Comenius University and CNRS/University of Nice-Sophia Antipolis}
\address[a]{Department of Applied Mathematics\\
\quad and Statistics\\
Faculty of Mathematics, Physics\\
\quad and Informatics\\
Comenius University\\
Bratislava\\
Slovakia\\
\printead{e1}}
\address[b]{Laboratoire I3S\\
CNRS/Universit\'{e} de Nice Sophia-Antipolis\\
B\^at. Euclide, Les Algorithmes, BP 121\\
2000 route des Lucioles\\
06903 Sophia Antipolis cedex\\
France\\
\printead{e2}}
\end{aug}

\received{\smonth{10} \syear{2013}}
\revised{\smonth{3} \syear{2014}}


\begin{abstract}
Among the major difficulties that one may encounter when estimating
parameters in a nonlinear regression model are the nonuniqueness of
the estimator, its instability with respect to small perturbations of
the observations and the presence of local optimizers of the estimation
criterion.

We show that these estimability issues can be taken into account at the
design stage, through the definition of suitable design criteria.
Extensions of $E$-, $c$- and $G$-optimality criteria are considered,
which when evaluated at a given $\mt^0$ (local optimal design), account
for the behavior of the model response $\eta(\mt)$ for $\mt$ far from
$\mt^0$. In particular, they ensure some protection against
close-to-overlapping situations where $\|\eta(\mt)-\eta(\mt^0)\|$ is
small for some $\mt$ far from $\mt^0$. These extended criteria are
concave and necessary and sufficient conditions for optimality
(equivalence theorems) can be formulated. They are not differentiable, but
when the design space is finite and the set $\Theta$ of admissible
$\mt
$ is discretized, optimal design forms a linear programming problem
which can be solved directly or via relaxation when $\Theta$ is just
compact. Several examples are presented.
\end{abstract}

%

\end{frontmatter}

\section{Introduction}

We consider a nonlinear regression model with observations
\[
y_i = y(x_i)=\eta(x_i,\mtb)+
\mve_i,\qquad  i=1,\ldots,N,
\]
where the errors $\mve_i$ satisfy $\Ex(\mve_i)=0$, $\operatorname
{var}(\mve_i)=\ms^2$
and $\cov(\mve_i,\mve_j)=0$ for $i\neq j$, $i,j=1,\ldots,N$, and the
true value $\mtb$ of the vector of model parameter $\mt$ belongs to~$\Theta$,
a compact subset of $\mathbb{R}^p$ such that $\Theta
\subset
\overline{\inter(\Theta)}$, the closure of the interior of~$\Theta
$. In
a vector notation, we write
%
\begin{equation}
\label{regression-model-vector}
 \yb=\eta_X(\mtb)+\mve \qquad\mbox{with } \Ex(\mve)=\0b,
\Var(\mve )=\ms^2 \Ib_N,
\end{equation}
where $\eta_X(\mt)=(\eta(x_1,\mt),\ldots,\eta(x_N,\mt))^\top $,
$\yb
=(y_1,\ldots,y_N)^\top $, $\mve=(\mve_1,\ldots,\mve_N)^\top $,
and $X$
denotes the $N$-point exact design $(x_1,\ldots,x_N)$.
The more general nonstationary (heteroscedastic) case where
$\operatorname{var}(\mve
_i)=\ms^2(x_i)$ can easily be transformed into the model (\ref
{regression-model-vector}) with $\ms^2=1$ via the division of
$y_i$ and $\eta(x_i,\mt)$ by $\ms(x_i)$.
We suppose that $\eta(x,\mt)$ is twice continuously differentiable with
respect to $\mt\in\inter(\Theta)$ for any $x\in\SX$, a compact subset
of $\mathbb{R}^d$. The model is assumed to be identifiable over $\SX$;
that is, we suppose that
%
\begin{equation}
\label{identifiable} \eta\bigl(x,\mt'\bigr)=\eta(x,\mt) \qquad
\mbox{for all } x
\in\SX \Longrightarrow\mt '=\mt.
\end{equation}

We shall denote by $\Xi$ the set of design measures $\xi$, that is, of
probability measures on $\SX$. The information matrix (for $\ms^2=1$)
for the design $X$ at $\mt$ is
\[
\Mb(X,\mt)= \sum_{i=1}^N
\frac{\mp\eta(x_i,\mt)}{\mp\mt} \frac
{\mp\eta
(x_i,\mt)}{\mp\mt^\top }
\]
and, for any $\xi\in\Xi$, we shall write
\[
\Mb(\xi,\mt)=\int_\SX\bigl[\mp\eta(x,\mt)/\mp\mt\bigr]
\bigl[\mp\eta(x,\mt )/\mp\mt^\top \bigr] \xi(\D x).
\]
%
Denoting $\xi_N=(1/N)\sum_{i=1}^N \delta_{x_i}$ the empirical design
measure associated with $X$, with $\delta_x$ the delta measure at $x$,
we have $\Mb(X,\mt)=N \Mb(\xi_N,\mt)$. Note that \eqref{identifiable}
implies the existence of a $\xi\in\Xi$ satisfying the Least-Squares
(LS) estimability condition
%
\begin{equation}
\label{estimability} \eta\bigl(x,\mt'\bigr)=\eta(x,\mt) \qquad\xi\mbox{-almost
everywhere} \Longrightarrow \mt'=\mt.
\end{equation}

Given an exact $N$-point design $X$, the set of all hypothetical means
of the observed vectors $\yb$ in the sample space $\mathbb{R}^N$ forms
the expectation surface
$\SSS_\eta=\{\eta_X(\mt)\dvtx\mt\in\Theta\}$.
Since $\eta_X(\mt)$ is supposed to have continuous first- and
second-order derivatives in $\inter(\Theta)$, $\SSS_\eta$ is a smooth
surface in $\mathbb{R}^N$ with a (local) dimension given by $r=\rank
[\mp
\eta_X(\mt)/\mp\mt^\top ]$. If $r=p$ (which means full rank),
the model~(\ref{regression-model-vector}) is said regular.
In regular models with no overlapping of $\SSS_\eta$, that is, when
$\eta_X(\mt) = \eta_X(\mt')$ implies $\mt= \mt'$, the LS estimator
%
\begin{equation}
\label{theta_LS} \mth_{\mathrm{LS}} = \mth_{\mathrm{LS}}^N = \arg\min
_{\mt\in\Theta} \bigl\|\yb-\eta _X(\mt)\bigr\|^2
\end{equation}
is uniquely defined with probability one (w.p.1). Indeed, when the
distributions of errors $\mve_i$ have probability densities (in the
standard sense) it can be proven that $\eta[\mth_{\mathrm{LS}}(\yb)]$ is unique
w.p.1; see \citet{Pazman84b} and \citeauthor{Pazman93}
(\citeyear{Pazman93}), page~107. However,
there is still a positive probability that the function $\mt
\longrightarrow\|\yb-\eta_X(\mt)\|^2$ has a local minimizer different
from the global one when the regression model is intrinsically curved
in the sense of \citet{BatesW80}, that is, when $\SSS_\eta$ is a curved
surface in $\mathbb{R}^N$; see \citeauthor{Demidenko00}
(\citeyear{Demidenko89,Demidenko00}).
Moreover, a curved surface can ``almost overlap''; that is, there may
exist points $\mt$ and $\mt'$ in $\Theta$ such that $\|\mt'-\mt\|$ is
large but $\|\eta_X(\mt')-\eta_X(\mt)\|$ is small (or even equals zero
in case of strict overlapping). This phenomenon can cause serious
difficulties in parameter estimation, leading to instabilities of the
estimator, and one should thus attempt to reduce its effects by
choosing an adequate experimental design.
Classically, those issues are ignored at the design stage and the
experiment is chosen on the basis of asymptotic local properties of the
estimator.
Even when the design relies on small-sample properties of the
estimator, like in \citet{PazmanPa92,GauchiP2006}, a nonoverlapping
assumption is used [see \citet{Pazman93}, pages~66 and 157] which permits
to avoid the aforementioned difficulties. Note that putting
restrictions on curvature measures is not enough: consider the case
$\dim(\mt)=1$ with the overlapping $\SSS_\eta$ formed by a circle of
arbitrarily large radius, and thus arbitrarily small curvature (see the
example in Section~\ref{S:motivating-ex} below).

Important and precise results are available concerning the construction
of subsets of $\Theta$ where such difficulties are guaranteed not to
occur; see, for example,
\citeauthor{Chavent83}
(\citeyear{Chavent83,Chavent90,Chavent91});
however, their exploitation for choosing adequate designs is far from
straightforward. Also, the construction of designs with restricted
curvatures, as proposed by \citet{ClydeC2002}, is based on the
curvature measures of \citet{BatesW80} and uses derivatives of $\eta
_X(\mt)$ at a certain $\mt$; this local approach is unable to catch the
problem of overlapping for two points that are distant in the parameter
space. Other design criteria using a second-order development of the
model response, or an approximation of the density of $\mth_{\mathrm{LS}}$
[\citet{HamiltonW85,PPa94}], are also inadequate.

The aim of this paper is to present new optimality criteria for optimum
design in nonlinear regression models that may reduce such effects,
especially overlapping, and are at the same time closely related to
classical optimality criteria like $E$-, $c$- or $G$-optimality (in
fact, they coincide with those criteria when the regression model is
linear). Classical optimality criteria focus on efficiency, that is,
aim at ensuring a precise estimation of $\mt$, asymptotically, provided
that the model is locally identifiable at $\mt$. On the other hand, the
new extended criteria account for the global behavior of the model and
enforce identifiability.

An elementary example is given in the next section and illustrates the
motivation of our work. The criterion of extended $E$-optimality is
considered in Section~\ref{S:extended-Eopt}; its main properties are
detailed and algorithms for the construction of optimal designs are
presented. Sections~\ref{S:E-c} and \ref{S:E-G} are, respectively,
devoted to the criteria of extended $c$-optimality and extended
$G$-optimality. Several illustrative examples are presented in
Section~\ref{S:example}. Section~\ref{S:Extensions} suggests some
extensions and further developments and Section~\ref{S:conclusions} concludes.

\section{An elementary motivating example}\label{S:motivating-ex}

\begin{example}\label{ex1}
Suppose that $\mt\in\Theta=[0,1]$ and that, for any design point
$x=(t, u)^\top \in\SX=\{0,\pi/2\}\times[0,u_{\max}]$, we have
\[
\eta(x,\mt)=r \cos(t-u \mt),
\]
with $r$ a known positive constant. We take $u_{\max}=7 \pi/4$; the
difficulties mentioned below are even more pronounced for values of
$u_{\max}$ closer to $2\pi$. We shall consider exclusively two-point
designs $X=(x_1,x_2)$ of the form
\[
x_1=(0, u)^\top,\qquad  x_2=(\pi/2,
u)^\top
\]
and denote $\nu_u$ the associated design measure, $\nu_u=(1/2)[\delta
_{x_1}+\delta_{x_2}]$. We shall look for an optimal design, that is, an
optimal choice of $u\in[0,u_{\max}]$, where optimality is considered in
terms of information.
\end{example}

\begin{figure}

\includegraphics{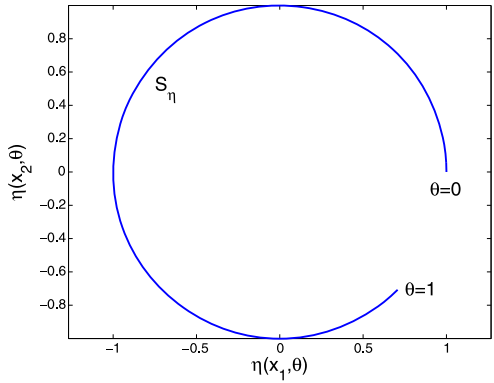}

\caption{Expectation surface $\SSS_\eta$ for $\mt\in\Theta=[0,1]$,
$r=1$ and $u=u_{\max}=7 \pi/4$.}
\label{F:S_eta}
\end{figure}

It is easy to see that for any design $\nu_u$ we have
\[
\eta_X(\mt)=\pmatrix{
\eta(x_1,\mt)
\cr
\eta(x_2,\mt)}
 = \pmatrix{
r \cos(u \mt)
\cr
r \sin(u \mt)}.
\]
The expectation surface is then an arc of a circle, with central angle
$u$; see Figure~\ref{F:S_eta} for the case $u=u_{\max}=7\pi/4$. The
model is nonlinear but parametrically linear since the information
matrix $M(X,\mt)$ for $\ms^2=1$ (here scalar since $\mt$ is scalar)
equals $r^2 u^2$ and does not depend on $\mt$. Also, the intrinsic
curvature (see Section~\ref{S:example}) is constant and equals $1/r$, and
the model is also almost intrinsically linear if $r$ gets large.

Any classical optimality criterion ($A$-, $D$-, $E$-) indicates that
one should observe at $u=u_{\max}$, and setting a constraint on the
intrinsic curvature is not possible here. However, if the true value of
$\mt$ is $\mtb=0$ and $\ms^2$ is large enough, there is a chance that
the LS estimator will be $\mth_{\mathrm{LS}}=1$, and thus very far from $\mtb$;
see Figure~\ref{F:S_eta}. The situation gets even worse if $u_{\max}$
gets closer to $2 \pi$, since $\SSS_\eta$ then almost overlaps.\vadjust{\goodbreak}

Now, consider $H_E(\nu_u,\mt)=(1/2) \|\eta_X(\mt)-\eta_X(\mt^0)\|
^2/|\mt
-\mt^0|^{2}$, see \eqref{H_E}, with $\mt^0=0$.
For all $u\in[0,u_{\max}]$, the minimum of $H_E(\nu_u,\mt)$ with
respect to $\mt\in\Theta$ is obtained at $\mt=1$, $H_E(\nu_u,1)=r^2
[1-\cos(u)]$ is then maximum in $[0,u_{\max}]$ for $u=u_*=\pi$. This
choice $u=u_*$ seems preferable to $u=u_{\max}$ since the expectation
surface $\SSS_\eta$ is then a half-circle, so that $\eta_X(0)$ and
$\eta
_X(1)$ are as far away as possible. On the other hand, as shown in
Section~\ref{S:extended-Eopt}, $\min_{\mt\in\Theta} H_E(\nu_u,\mt)$
possesses most of the attractive properties of classical optimality
criteria and even coincides with one of them in linear models.

\begin{figure}

\includegraphics{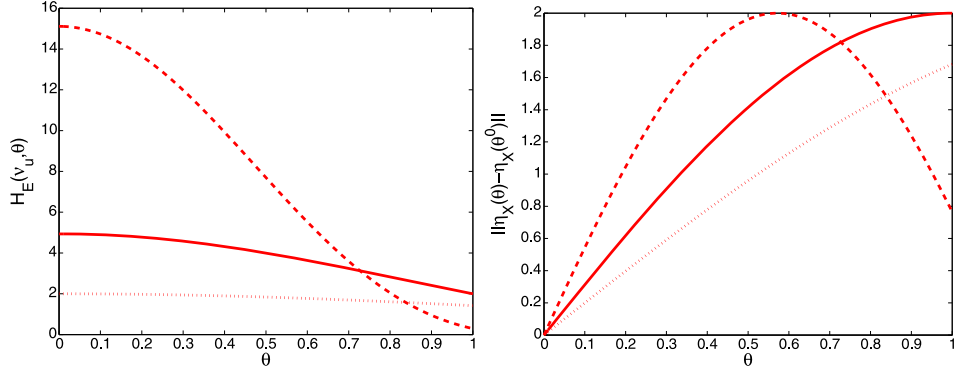}

\caption{$H_E(\nu_u,\mt)$ (\emph{left}) and $\|\eta_X(\mt)-\eta_X(\mt
^0)\|$
(\emph{right}) as functions of
$\mt\in\Theta=[0,1]$ for $r=1$, $u=2$ (\emph{dotted line}),
$u=u_{\max
}=7\pi/4$ (\emph{dashed line}) and $u=u_*=\pi$ (\emph{solid line}).}
\label{F:H_x_theta}
\end{figure}

Figure~\ref{F:H_x_theta}-left shows $H_E(\nu_u,\mt)$ as a function of
$\mt$ for three values of $u$ and illustrates the fact that the minimum
of $H_E(\nu_u,\mt)$ with respect to $\mt\in\Theta$ is maximized for
$u=u_*$. Figure~\ref{F:H_x_theta}-right shows that the design with
$u=u_{\max}$ (\emph{dashed line}) is optimal locally at $\mt=\mt
^0$, in
the sense that it yields the fastest increase of $\|\eta_X(\mt)-\eta
_X(\mt^0)\|$ as $\mt$ slightly deviates from $\mt^0$. On the other
hand, $u=\pi$ maximizes $\min_{\mt\in\Theta} \|\eta_X(\mt)-\eta
_X(\mt
^0)\|/|\mt-\mt^0|$ (\emph{solid line}) and realizes a better protection
against the folding effect of $\SSS_\eta$, at the price of a slightly
less informative experiment for $\mt$ close to $\mt^0$. Smaller values
of $u$ (\emph{dotted line}) are worse than~$u_*$, both locally for
$\mt
$ close to $\mt^0$ and globally in terms of the folding of $\SSS_\eta$.

The rest of the paper will formalize these ideas and show how to
implement them for general nonlinear models through the definition of
suitable design criteria that can be easily optimized.

\section{Extended (globalized) $E$-optimality}\label{S:extended-Eopt}

\subsection{Definition of \texorpdfstring{$\phi_{eE}(\cdot)$}{phi eE(.)}}\label{S:eE-def}

Take a fixed point $\mt^0$ in $\Theta$ and denote
%
\begin{equation}
\label{H_E} H_E(\xi,\mt)= H_E\bigl(\xi,\mt;
\mt^0\bigr) = \frac{\|\eta(\cdot,\mt
)-\eta(\cdot,\mt^0)\|_\xi^2}{\|\mt-\mt^0\|^{2}},\vadjust{\goodbreak}
\end{equation}
where $\llVert \cdot\rrVert _\xi$ denotes the norm in $\SL_2(\xi
)$; that
is, $\llVert l \rrVert _\xi= [\int_{\SX} l^2(x) \xi(\D x)
]^{1/2}$ for any $l \in\SL_2(\xi)$. When $\xi$ is a discrete measure,
like in the examples considered in the paper, then $\llVert l\rrVert
_\xi
^2$ is simply the sum $\sum_{x: \xi(\{x\})>0} \xi(\{x\}) l^2(x)$.

The extended $E$-optimality criterion is defined by
%
\begin{equation}
\label{phi_eE} \phi_{eE}(\xi)=\phi_{eE}\bigl(\xi;
\mt^0\bigr) = \min_{\mt\in\Theta} H_E(\xi,\mt),
\end{equation}
to be maximized with respect to the design measure $\xi$.

In a nonlinear regression model $\phi_{eE}(\cdot)$ depends on the value
chosen for $\mt^0$ and can thus be considered as a local optimality
criterion.\vadjust{\goodbreak} On the other hand, the criterion is global in the sense that
it depends on the behavior of $\eta(\cdot,\mt)$ for $\mt$ far from
$\mt
^0$. This (limited) locality can be removed by considering $\phi
_{MeE}(\xi)=\min_{\mt^0\in\Theta}\phi_{eE}(\xi;\mt^0)$ instead
of \eqref
{phi_eE}, but only the case of $\phi_{eE}(\cdot)$ will be detailed in
the paper, the developments being similar for $\phi_{MeE}(\cdot)$; see
Section~\ref{S:Maximin-extended-Eopt}.

For a linear regression model with
$\eta(x,\mt)=\mathbf{f}^\top (x)\mt+v(x)$
and $\Theta=\mathbb{R}^p$, for any $\mt^0$ and any $\xi\in\Xi$,
we have
$\|\eta(\cdot,\mt)-\eta(\cdot,\mt^0)\|_\xi^2 = (\mt-\mt
^0)^\top
\Mb(\xi
)(\mt-\mt^0)$, so that
\[
\phi_{eE}(\xi) = \min_{\mt-\mt^0\in\mathbb{R}^p} \frac{(\mt
-\mt^0)^\top \Mb
(\xi)(\mt-\mt^0)}{\|\mt-\mt^0\|^2} =
\ml_{\min}\bigl[\Mb(\xi)\bigr],
\]
the minimum eigenvalue of $\Mb(\xi)$, and corresponds to the
$E$-optimality criterion.

For a nonlinear model with $\Theta=\SB(\mt^0,\rho)$, the ball with
center $\mt^0$ and radius $\rho$, direct calculation shows that
%
\begin{equation}
\label{eE_rho-small} \lim_{\rho\ra0} \phi_{eE}\bigl(\xi;
\mt^0\bigr) = \ml_{\min}\bigl[\Mb\bigl(\xi,\mt ^0
\bigr)\bigr].
\end{equation}
In a nonlinear regression model with larger $\Theta$, the determination
of an optimum design $\xi_{eE}^*$ maximizing $\phi_{eE}(\xi)$ ensures
some protection against $\|\eta(\cdot,\mt)-\eta(\cdot,\mt^0)\|
_\xi$
being small\vspace*{1pt} for some $\mt$ far from $\mt^0$. In particular, when $\mt
^0\in\inter(\Theta)$ then $\phi_{eE}(\xi;\mt^0)=0$ if either $\Mb
(\xi,\mt^0)$ is singular or $\|\eta(\cdot,\mt)-\eta(\cdot,\mt^0)\|
_\xi=0$
for some $\mt\neq\mt^0$. Therefore, under the condition \eqref
{identifiable}, $\xi_{eE}^*$ satisfies the estimability condition~\eqref
{estimability} at $\mt=\mt^0$
and is necessarily nondegenerate, that is, $\Mb(\xi_{eE}^*,\mt^0)$ is
nonsingular, when $\mt^0\in\inter(\Theta)$ (provided that there exists
a nondegenerate design in $\Xi$). Notice that \eqref{eE_rho-small}
implies that $\phi_{eE}(\xi;\mt^0)\leq\ml_{\min}[\Mb(\xi,\mt
^0)]$ when
$\Theta$ contains some open neighborhood of $\mt^0$. In contrast with
the $E$-optimality criterion, maximizing $\phi_{eE}(\xi;\mt^0)$ in
nonlinear models does not require computation of the derivatives of
$\eta(x,\mt)$ with respect to $\mt$ at $\mt^0$; see the algorithms
proposed in Sections~\ref{S:eE-LP} and~\ref{S:eE-relax}. Also note that
the influence of points that are very far from $\mt^0$ can be
suppressed by modification of the denominator of \eqref{H_E} without
changing the relation with $E$-optimality; see Section~\ref{S:extra K}.

Before investigating properties of $\phi_{eE}(\cdot)$ as a criterion
function for optimum design in the next section, we state a property
relating $\phi_{eE}(\xi)$ to the localization of the LS estimator
$\mth_{\mathrm{LS}}$.

\begin{theo}\label{Th:localization-thLS}
For any given $\mt\in\Theta$, the LS estimator $\mth_{\mathrm{LS}}$ given by
\eqref{theta_LS} in the model \eqref{regression-model-vector} satisfies
\[
\mth_{\mathrm{LS}} \in\Theta\cap\SB \biggl(\mt,\frac{2 \|\yb-\eta_X(\mt
)\|}{\sqrt
{N} \sqrt{\phi_{eE}(\xi_N;\mt)}} \biggr),
\]
with $\xi_N$ the empirical measure associated with the design $X$ used
to observe~$\yb$.
\end{theo}

\begin{pf} The result follows from the following chain of
inequalities:
%
\begin{eqnarray}
\label{bound-mthLS} \|\mth_{\mathrm{LS}}-\mt\| &\leq& \frac{\|\eta(\cdot,\mth_{\mathrm{LS}})-\eta
(\cdot,\mt)\|
_{\xi_N}}{\sqrt{\phi_{eE}(\xi_N;\mt)}} =
\frac{\|\eta_X(\mth
_{\mathrm{LS}})-\eta
_X(\mt)\|}{\sqrt{N} \sqrt{\phi_{eE}(\xi_N;\mt)}}
\nonumber
\\
& \leq&\frac{\|\yb-\eta_X(\mth_{\mathrm{LS}})\|+\|\yb-\eta_X(\mt)\|
}{\sqrt{N}
\sqrt{\phi_{eE}(\xi_N;\mt)}} \leq\frac{2 \|\yb-\eta_X(\mt)\|
}{\sqrt{N}
\sqrt{\phi_{eE}(\xi_N;\mt)}}.
\end{eqnarray}
\upqed\end{pf}

Note that although the bound \eqref{bound-mthLS} is tight in
general nonlinear situations (due to the possibility that $\SSS_\eta$
overlaps), it is often pessimistic. In particular, in the linear
regression model
$\eta(x,\mt)=\mathbf{f}^\top (x)\mt+v(x)$,
direct calculation gives
\[
\|\mth_{\mathrm{LS}}-\mt\| \leq\sqrt{\ml_{\max}\bigl[\bigl(
\Fb^\top \Fb\bigr)^{-1}\bigr]} \bigl\| \yb-\eta _X(
\mt)\bigr\| = \frac{\|\yb-\eta_X(\mt)\|}{\sqrt{N} \sqrt{\phi
_{eE}(\xi
_N)}},
\]
where $\Fb$ is the $N\times p$ matrix with $i$th line equal to
$\mathbf{f}^\top (x_i)$.
We also have $\|\mth_{\mathrm{LS}}-\mt\| \leq\|\yb-\eta_X(\mt)\| /[\sqrt{N}
\sqrt{\phi_{eE}(\xi_N,\mt)}]$ in
intrinsically linear models (with a flat expectation surface $\SSS
_\eta
$) since then $\|\eta_X(\mth_{\mathrm{LS}})-\eta_X(\mt)\|\leq\|\yb-\eta
_X(\mt)\|$.

In the following, we shall omit the dependence in $\mt^0$ and simply
write $\phi_{eE}(\xi)$ for $\phi_{eE}(\xi;\mt^0)$ when there is no
ambiguity.

\subsection{Properties of \texorpdfstring{$\phi_{eE}(\cdot)$}{phi eE(.)}}\label{S:eE-prop}
As the minimum of linear functions of $\xi$, $\phi_{eE}(\cdot)$ is
\emph
{concave}: for all $\xi,\nu\in\Xi$ and all $\ma\in[0,1]$, $\phi
_{eE}[(1-\ma)\xi+\ma\nu]\geq(1-\ma)\phi_{eE}(\xi)+\ma\phi
_{eE}(\nu)$.
It is also \emph{positively homogeneous}: $\phi_{eE}(a\xi)=a \phi
_{eE}(\xi)$ for all $\xi\in\Xi$ and $a>0$; see, for example, \citet{Pukelsheim93}, Chapter~5.
The criterion of ${eE}$-efficiency can then be
defined as
\[
\SE_{eE}(\xi)= \frac{\phi_{eE}(\xi)}{\phi_{eE}(\xi_{eE}^*)},\qquad \xi\in \Xi,
\]
where $\xi_{eE}^*$ maximizes $\phi_{eE}(\xi)$.

The concavity of $\phi_{eE}(\cdot)$ implies the existence of
directional derivatives and, due to the linearity in $\xi$ of $H_E(\xi,\mt)$, we
have the following; see, for example, \citet{DemyanovM74}.\vadjust{\goodbreak}

\begin{theo}
For any $\xi,\nu\in\Xi$, the directional derivative of the criterion
$\phi_{eE}(\cdot)$ at $\xi$ in the direction $\nu$ is given by
\[
F_{\phi_{eE}}(\xi;\nu)=\min_{\mt\in\Theta_E(\xi)} H_E(\nu,\mt
) - \phi _{eE}(\xi),
\]
where
$\Theta_E(\xi) =  \{\mt\in\Theta\dvtx H_E(\xi,\mt) = \phi
_{eE}(\xi)
 \}$.
\end{theo}

Note that we can write
$F_{\phi_{eE}}(\xi;\nu)=\min_{\mt\in\Theta_E(\xi)} \int_\SX
\Psi
_{eE}(x,\mt,\xi) \nu(\D x)$,
where
%
\begin{equation}
\label{Psi_Est,mt} \Psi_{eE}(x,\mt,\xi) = \frac{[\eta(x,\mt)-\eta(x,\mt^0)]^2 - \|
\eta
(\cdot,\mt)-\eta(\cdot,\mt^0)\|_\xi^2}{\|\mt-\mt^0\|^2}.
\end{equation}
Due to the concavity of $\phi_{eE}(\cdot)$, a necessary and sufficient
condition for the optimality of a design measure $\xi_{eE}^*$ is that
%
\begin{equation}
\label{CNS_eE1} \sup_{\nu\in\Xi} F_{\phi_{eE}}\bigl(
\xi_{eE}^*;\nu\bigr) = 0,
\end{equation}
a condition often called ``equivalence theorem'' in optimal design
theory; see, for example, \citet{Fedorov72,Silvey80}.
An equivalent condition is as follows.

\begin{theo}\label{Th:Eq-eE}
A design $\xi_{eE}^*\in\Xi$ is optimal for $\phi_{eE}(\cdot)$ if and
only if
%
\begin{eqnarray}
\label{CNS_eE2} \max_{x\in\SX} \int_{\Theta_E(\xi_{eE}^*)}
\Psi_{eE}(x,\mt,\xi ) \mu ^*(\D\mt) = 0
\nonumber
\\[-8pt]
\\[-8pt]
\eqntext{\mbox{for some measure $\mu^*
\in\SM\bigl[\Theta_E\bigl(\xi _{eE}^*\bigr)\bigr]$},}
\end{eqnarray}
the set of probability measures on $\Theta_E(\xi_{eE}^*)$.
\end{theo}

\begin{pf}
This is a classical result for maximin design problems; see, for
example, \citet{FedorovH97}, Section~2.6. We have
%
\begin{eqnarray}
\label{*minimax}
0 &\leq&\sup_{\nu\in\Xi} F_{\phi_{eE}}(\xi;\nu) \nonumber\\
&=&
\sup_{\nu
\in\Xi} \min_{\mt\in\Theta_E(\xi)} \int
_\SX\Psi_{eE}(x,\mt,\xi) \nu (\D x)
\nonumber
\\
&=& \sup_{\nu\in\Xi} \min_{\mu\in\SM[\Theta_E(\xi)]} \int
_\SX\int_{\Theta_E(\xi)} \Psi_{eE}(x,\mt,
\xi) \mu(\D\mt) \nu(\D x)\nonumber
\\
&=& \min_{\mu\in\SM[\Theta_E(\xi)]} \sup_{\nu\in\Xi} \int
_\SX\int_{\Theta_E(\xi)} \Psi_{eE}(x,\mt,
\xi) \mu(\D\mt) \nu(\D x)
\nonumber
\\
&= &\min_{\mu\in\SM[\Theta_E(\xi)]} \max_{x\in\SX} \int
_{\Theta_E(\xi
)} \Psi_{eE}(x,\mt,\xi) \mu(\D\mt).
\end{eqnarray}
Therefore, the necessary and sufficient condition (\ref{CNS_eE1}) can
be written as (\ref{CNS_eE2}).\vadjust{\goodbreak}~%
\end{pf}

One should notice that $\sup_{\nu\in\Xi} F_{\phi_{eE}}(\xi;\nu)$ is
generally not obtained for $\nu$ equal to a one-point (delta) measure,
which prohibits the usage of classical vertex-direction algorithms for
optimizing $\phi_{eE}(\cdot)$. Indeed, the minimax problem (\ref
{*minimax}) has generally several solutions $x^{(i)}$ for $x$,
$i=1,\ldots,s$, and the optimal $\nu^*$ is then a linear combination
$\sum_{i=1}^s w_i \delta_{x^{(i)}}$, with $w_i\geq0$ and $\sum_{i=1}^s
w_i=1$; see \citet{PHWc91} for developments on a similar difficulty in
$T$-optimum design for model discrimination. This property, due to the
fact that $\phi_{eE}(\cdot)$ is not differentiable, has the important
consequence that the determination of a maximin-optimal design cannot
be obtained via standard design algorithms used for differentiable criteria.

To avoid that difficulty, a regularized version $\phi_{eE,\ml}(\cdot)$
of $\phi_{eE}(\cdot)$ is considered in
\citet{PP2013}, Sections~7.7.3 and 8.3.2,
with the property that $\lim_{\ml\ra\infty}\phi
_{eE,\ml
}(\xi)=\phi_{eE}(\xi)$ for any $\xi\in\Xi$ (the convergence being
uniform when $\Theta$ is a finite set), $\phi_{eE}(\cdot)$ is concave
and such that $\sup_{\nu\in\Xi} F_{\phi_{eE,\ml}}(\xi;\nu)$ is obtained
when $\nu$ is the delta measure $\delta_{x^*}$ at some $x^*\in\SX$
(depending on~$\xi$). However, although $\phi_{eE,\ml}(\cdot)$ is
smooth for any finite $\ml$, its maximization tends to be badly
conditioned for large $\ml$.

In the next section, we show that optimal design for $\phi_{eE}(\cdot)$
reduces to linear programming when $\Theta$ and $\SX$ are finite. This
is an important property. An algorithm based on a relaxation of the
maximin problem is then considered in Section~\ref{S:eE-relax} for the
case where $\Theta$ is compact.

\subsection{Optimal design via linear-programming (\texorpdfstring{$\Theta$}{Theta} is
finite)}\label{S:eE-LP}
To simplify the construction of an optimal design, one may take $\Theta
$ as a finite set, $\Theta=\Theta^{(m)}=\{\mt^{(1)},\mt
^{(2)},\ldots,\mt
^{(m)}\}$; $\phi_{eE}(\xi)$ can then be written as $\phi_{eE}(\xi
)= \min_{j=1,\ldots,m} H_E(\xi,\mt^{(j)})$, with $H_E(\xi,\mt)$ given by
(\ref
{H_E}). If the design space $\SX$ is also finite, with $\SX=\{
x^{(1)},x^{(2)},\ldots,x^{(\ell)}\}$, then the determination of an
optimal design measure for $\phi_{eE}(\cdot)$ amounts to the
determination of a scalar $t$ and of a vector of weights $\wb
=(w_1,w_2,\ldots,w_\ell)^\top $, $w_i$ being allocated at $x^{(i)}$ for
each $i=1,\ldots,\ell$, such that $\cb^\top [\wb^\top, t]^\top
$ is
maximized, with $\cb=(0,0,\ldots,0,1)^\top $ and $\wb$ and $t$ satisfying
the constraints
%
\begin{eqnarray}
\label{CTLP}  \sum_{i=1}^\ell
w_i&=&1,\qquad  w_i \geq0, i=1,\ldots,\ell,
\nonumber
\\[-8pt]
\\[-8pt]
\nonumber
\sum_{i=1}^\ell w_i
h_i\bigl(\mt^{(j)}\bigr) &\geq& t,\qquad j=1,\ldots,m,
\end{eqnarray}
where we denoted
%
\begin{equation}
\label{h_i} h_i(\mt) = \frac{[\eta(x^{(i)},\mt)-\eta(x^{(i)},\mt^0)]^2}{\|\mt
-\mt
^0\|^{2}}.
\end{equation}
This is a linear programming (LP) problem, which can easily be solved
using standard methods (for instance, the simplex algorithm), even for
large $m$ and $\ell$. We shall denote by $(\hat\wb, \hat
t)=\mathrm{LP}_{eE}(\SX,\Theta^{(m)})$ the solution of this problem.

We show below how a compact subset $\Theta$ of $\mathbb{R}^p$ with
nonempty interior can be replaced by a suitable discretized version
$\Theta^{(m)}$ that can be enlarged iteratively.

\subsection{Optimal design via relaxation and the cutting-plane method
(\texorpdfstring{$\Theta$}{Theta} is a compact subset of $\mathbb{R}^p$)}\label{S:eE-relax}

Suppose now that $\SX$ is finite and that $\Theta$ is a compact subset
of $\mathbb{R}^p$ with nonempty interior. In the LP formulation above,
$(\wb, t)$ must satisfy an infinite number of constraints:
$\sum_{i=1}^\ell w_i h_i(\mt) \geq t$ for all $\mt\in\Theta$; see \eqref
{CTLP}. One may then use the method of \citet{ShimizuA80} and consider
the solution of a series of relaxed LP problems, using at step $k$ a
finite set of constraints only, that is, consider $\mt\in\Theta^{(k)}$
finite. Once a solution $(\wb^k,t^k)=\mathrm{LP}_{eE}(\SX,\Theta^{(k)})$ of this
problem is obtained, using a standard LP solver, the set $\Theta^{(k)}$
is enlarged to $\Theta^{(k+1)}=\Theta^{(k)}\cup\{\mt^{(k+1)}\}$ with
$\mt^{(k+1)}$ given by the constraint \eqref{CTLP} most violated by
$\wb
^k$, that is,
%
\begin{equation}
\label{theta_k+1} \mt^{(k+1)}=\arg\min_{\mt\in\Theta}
H_E\bigl(\wb^k,\mt\bigr),
\end{equation}
where with a slight abuse of notation, we write $H_E(\wb,\mt)=H_E(\xi,\mt)$; see \eqref{H_E}, when $\xi$ allocates mass $w_i$ at the support
point $x^{(i)}\in\SX$ for all $i$.
This yields the following algorithm for the maximization of $\phi
_{eE}(\cdot)$.

\begin{longlist}[(1)]
\item[(0)] Take any vector $\wb^0$ of nonnegative weights summing to
one, choose $\me>0$, set $\Theta^{(0)}=\varnothing$ and $k=0$.
\item[(1)] Compute $\mt^{(k+1)}$ given by (\ref{theta_k+1}), set
$\Theta
^{(k+1)}=\Theta^{(k)}\cup\{\mt^{(k+1)}\}$.
\item[(2)] Use a LP solver to determine $(\wb^{k+1}, t^{k+1}
)=\mathrm{LP}_{eE}(\SX,\Theta^{(k+1)})$.
\item[(3)] If $\Delta_{k+1}=t^{k+1} - \phi_{eE}(\wb^{k+1}) < \me$, take
$\wb^{k+1}$ as an $\me$-optimal solution and stop; otherwise
$k\leftarrow k+1$, return to step 1.
\end{longlist}

The optimal value $\phi_{eE}^*=\max_{\xi\in\Xi} \phi_{eE}(\xi)$ satisfies
\[
\phi_{eE}\bigl(\wb^{k+1}\bigr) \leq\phi_{eE}^* \leq
t^{k+1}
\]
at every iteration, so that $\Delta_{k+1}$ of step 3 gives an upper
bound on the distance to the optimum in terms of criterion value.


The algorithm can be interpreted in terms of the cutting-plane method.
Indeed, from \eqref{H_E} and \eqref{h_i} we have $H_E(\wb,\mt
^{(j+1)})=\sum_{i=1}^\ell w_i h_i(\mt^{(j+1)})$ for any vector of
weights $\wb$. From the definition of $\mt^{(j+1)}$ in \eqref
{theta_k+1}, we obtain
\begin{eqnarray*}
\phi_{eE}(\wb) \leq H_E\bigl(\wb,\mt^{(j+1)}\bigr)
&=& H_E\bigl(\wb^j,\mt ^{(j+1)}\bigr) + \sum
_{i=1}^\ell h_i\bigl(
\mt^{(j+1)}\bigr) \bigl\{\wb-\wb^j\bigr\}_i
\\
&=& \phi_{eE}\bigl(\wb^j\bigr) + \sum
_{i=1}^\ell h_i\bigl(\mt^{(j+1)}
\bigr) \bigl\{\wb-\wb ^j\bigr\} _i,
\end{eqnarray*}
so that the vector with components $h_i(\mt^{(j+1)})$, $i=1,\ldots,\ell
$, forms a subgradient of $\phi_{eE}(\cdot)$ at $\wb^j$, which we
denote $\nabla\phi_{eE}(\wb^j)$ below [it is sometimes called
supergradient since $\phi_{eE}(\cdot)$ is concave].
%
Each of the constraints
\[
\sum_{i=1}^\ell w_i
h_i\bigl(\mt^{(j+1)}\bigr) \geq t,
\]
used in the LP problem of step 2, with $j=0,\ldots,k$, can be written as
\[
\nabla^\top \phi_{eE}\bigl(\wb^j\bigr)\wb=
\phi_{eE}\bigl(\wb^j\bigr)+\nabla^\top \phi
_{eE}\bigl(\wb ^j\bigr) \bigl(\wb-\wb^j\bigr)
\geq t.
\]
Therefore, $\wb^{k+1}$ determined at step 2 maximizes the
piecewise-linear approximation
\[
\min_{j=0,\ldots,k} \bigl\{\phi_{eE}\bigl(\wb^j
\bigr)+\nabla^\top \phi_{eE}\bigl(\wb ^j\bigr)
\bigl(\wb -\wb^j\bigr)\bigr\}
\]
of $\phi_{eE}(\wb)$ with respect to the vector of weights $\wb$, and
the algorithm corresponds to the cutting-plane method of \citet{Kelley60}.

The only difficult step in the algorithm corresponds to the
determination of $\mt^{(k+1)}$ in \eqref{theta_k+1} when $\Theta$ is a
compact set. We found that the following simple procedure is rather efficient.
Construct a finite grid, or a space-filling design, $\SG^0$ in $\Theta
$. Then, for $k=0,1,2,\ldots$
%
\begin{equation}
\label{theta_k+1_V2} \cases{ %
\phantom{ii}\mathrm{(i)} & \quad$\mbox{compute }
\mth^{k+1}=\displaystyle\arg\min_{\mt'\in\SG^k} H_E\bigl(
\wb^k,\mt '\bigr);$
\vspace*{2pt}\cr
\phantom{i}\mathrm{(ii)} &\quad  $\mbox{perform a local minimization of } H_E\bigl(
\wb^k,\mt\bigr)$
\vspace*{2pt}\cr
&\quad $\mbox{with respect to } \mt\in\Theta, \mbox{ initialized at } \mth
^{k+1};$
\vspace*{2pt}\cr
&\quad $\mbox{let } \mt^{(k+1)} \mbox{ denote the solution};$
\vspace*{2pt}\cr
\mathrm{(iii)} &\quad $\mbox{set } \SG^{k+1}=\SG^k\cup\bigl\{
\mt^{(k+1)}\bigr\}.$}
\end{equation}
%
The optimal value $\phi_{eE}(\xi^*_{eE})$ can then be approximated by
$H_{E}(\wb^{k+1},\mt^{(k+2)})$ when the algorithm stops (step 3).

The method of cutting planes is known to have sometimes rather poor
convergence properties; see, for example,
\citeauthor{BonnansGLS2006} (\citeyear{BonnansGLS2006}), Chapter~9,
\citeauthor{Nesterov2004} (\citeyear{Nesterov2004}), Section~3.3.2. A
significant improvement consists in restricting the search for $\wb
^{k+1}$ at step 2 to some neighborhood of the best solution obtained so
far, which forms the central idea of bundle methods; see \citet
{LemarechalNN95}, \citeauthor{BonnansGLS2006} (\citeyear{BonnansGLS2006}),
Chapters~9--10. In particular,
the level method of \citeauthor{Nesterov2004}
(\citeyear{Nesterov2004}), Section~3.3.3, adds a
quadratic-programming step to each iteration of the cutting planes
algorithm presented above; one may refer for instance to
\citeauthor{PP2013} (\citeyear{PP2013}),
Section~9.5.3, for an application of the level method to
design problems.
Notice that any linear constraint on $\wb$ can easily be taken into
account in addition to those in \eqref{CTLP}, so that the method
directly applies to optimal design with linear cost-constraints; see,
for example, \citeauthor{FedorovL2014} (\citeyear{FedorovL2014}),
Section~4.2.


\section{Extended (globalized) $c$-optimality}\label{S:E-c}

\subsection{Definition and properties} \label{S:ec-def&prop}

Consider the case where one wants to estimate a scalar function of $\mt
$, denoted by $g(\mt)$, possibly nonlinear. We assume that
\[
\cb=\cb(\mt)=\frac{\mp g(\mt)}{\mp\mt} \bigg|_{\mt=\mt^0} \neq \0b.
\]
Denote
%
\begin{equation}
\label{H_c} H_c(\xi,\mt)= H_c\bigl(\xi,\mt;
\mt^0\bigr) = \frac{\|\eta(\cdot,\mt
)-\eta(\cdot,\mt^0)\|_\xi^2}{| g(\mt)- g(\mt^0)|^{2}}
\end{equation}
and consider the design criterion defined by
%
\begin{equation}
\label{phi_ec} \phi_{ec}(\xi)= \min_{\mt\in\Theta}
H_c(\xi,\mt),
\end{equation}
to be maximized with respect to the design measure $\xi$.

When $\eta(x,\mt)$ and the scalar function $ g(\mt)$ are both linear in
$\mt$, with $ g(\mt)=\cb^\top \mt$, we get
\[
\phi_{ec}(\xi) = \min_{\mt\in\Theta, \cb^\top (\mt-\mt^0)\neq
0} \frac
{(\mt-\mt^0)^\top \Mb(\xi) (\mt-\mt^0)}{[\cb^\top (\mt-\mt^0)]^2}
\]
and, therefore, $\phi_{ec}(\xi) = [\cb^\top \Mb^-(\xi)\cb
]^{-1}$, using
the well-known formula\break  $\cb^\top \Mb^-\cb=\max_{\ma\neq0} (\cb
^\top
\ma
)^2/(\ma^\top \Mb\ma)$; cf. \citeauthor{Harville97} (\citeyear{Harville97}), equation~(10.4).
Also, for a nonlinear model with $\Theta=\SB(\mt^0,\rho)$ and a design
$\xi$ such that $\Mb(\xi,\mt^0)$ has full rank, one has
\[
\lim_{\rho\ra0} \phi_{ec}(\xi) = \bigl[\cb^\top
\Mb^{-1}\bigl(\xi,\mt ^0\bigr)\cb\bigr]^{-1},
\]
which justifies that we consider $\phi_{ec}(\xi)$ as an \emph{extended
$c$-optimality criterion}.
At the same time, in a nonlinear situation with larger $\Theta$ the
determination of an optimal design $\xi_{ec}^*$ maximizing $\phi
_{ec}(\xi)$ ensures some protection against $\|\eta(\cdot,\mt)-\eta
(\cdot,\mt^0)\|_\xi^2$ being small for some $\mt$ such that $ g(\mt)$
is significantly different from $ g(\mt^0)$. The condition \eqref
{identifiable} guarantees the existence of a $\xi\in\Xi$ such that
$\phi
_{ec}(\xi)>0$, and thus the LS estimability of $ g(\mt)$ at $\mt^0$ for
$\xi_{ec}^*$, that is,
\[
\eta(x,\mt)=\eta\bigl(x,\mt^0\bigr), \qquad\xi_{ec}^*
\mbox{-almost everywhere } \Longrightarrow g(\mt)= g\bigl(\mt^0\bigr);
\]
see \citeauthor{PP2013} (\citeyear{PP2013}), Section~7.4.4.
When $\Theta$ contains an open
neighborhood of $\mt^0$, then $\phi_{ec}(\xi)\leq[\cb^\top \Mb
^-(\xi,\mt
^0)\cb]^{-1}$.

Similarly to $\phi_{eE}(\cdot)$, the criterion $\phi_{ec}(\cdot)$ is
concave and positively homogeneous; its concavity implies the existence
of directional derivatives.

\begin{theo}
For any $\xi,\nu\in\Xi$, the directional derivative of the criterion
$\phi_{ec}(\cdot)$ at $\xi$ in the direction $\nu$ is given by
\[
F_{\phi_{ec}}(\xi;\nu)=\min_{\mt\in\Theta_c(\xi)} H_c(\nu,\mt
) - \phi _{ec}(\xi),
\]
where
$\Theta_c(\xi) =  \{\mt\in\Theta\dvtx H_c(\xi,\mt) = \phi
_{ec}(\xi)
 \}$.
\end{theo}

A necessary and sufficient condition for the optimality of $\xi^*$
maximizing $\phi_{ec}(\cdot)$ is that $\sup_{\nu\in\Xi} F_{\phi
_{ec}}(\xi^*;\nu) = 0$, which yields an equivalence theorem similar to
Theorem~\ref{Th:Eq-eE}.

When both $\Theta$ and $\SX$ are finite, an optimal design for $\phi
_{ec}(\cdot)$ is obtained by solving a LP problem. Compared with
Section~\ref{S:eE-LP}, we simply need to substitute $H_c$ for $H_E$ and
use $h_i(\mt) = [\eta(x^{(i)},\mt)-\eta(x^{(i)},\mt^0)]^2/| g(\mt)-
g(\mt^0)|^{2}$, $i=1,\ldots,\ell$, instead of \eqref{h_i}. Also, a
relaxation method similar to that in Section~\ref{S:eE-relax} can be
used when $\Theta$ is a compact subset of $\mathbb{R}^p$.

\section{Extended (globalized) $G$-optimality}\label{S:E-G}
Following the same lines as above, we can also define an extended
$G$-optimality criterion by
\[
\phi_{eG}(\xi) = \min_{\mt\in\Theta} \frac{\|\eta(\cdot,\mt
)-\eta(\cdot,\mt^0)\|_\xi^2}{\max_{x\in\SX} [\eta(x,\mt)-\eta(x,\mt
^0)
]^2}.
\]
The fact that it corresponds to the $G$-optimality criterion for a
linear model can easily be seen, noticing that in the model (\ref
{regression-model-vector}) with $\eta(x,\mt)=\mathbf{f}^\top (x)\mt+v(x)$
we have
\begin{eqnarray*}
\biggl\{\sup_{x\in\SX} \frac{N}{\ms^2} \operatorname{var}
\bigl[ \mathbf{f}^\top (x)\mth_{\mathrm{LS}} \bigr] \biggr
\}^{-1} &=& \inf_{x\in\SX} \bigl[\mathbf{f}^\top
(x) \Mb ^{-1}(\xi_N) \mathbf{f} (x)\bigr]^{-1}
\\
&=& \inf_{x\in\SX} \inf_{\ub\in\mathbb{R}^p, \ub^\top \mathbf{f}(x)
\neq0}
\frac
{\ub^\top \Mb(\xi_N)\ub}{[\mathbf{f}^\top (x)\ub]^2}
\\
&=& \inf_{\ub\in\mathbb{R}^p} \frac{\ub^\top \Mb(\xi_N)\ub
}{\max_{x\in\SX
}[\mathbf{f}^\top (x)\ub]^2},
\end{eqnarray*}
where $\xi_N$ denotes the empirical design measure corresponding to
$X$, assumed to be nonsingular, and the second equality follows from
\citeauthor{Harville97} (\citeyear{Harville97}), equation (10.4).
The equivalence theorem of \citet{KieferW60} indicates that $D$- and
$G$-optimal designs coincide; therefore, $D$-optimal designs are
optimal for $\phi_{eG}(\cdot)$ in linear models. Moreover, the optimum
(maximum) value of
$\phi_{eG}(\xi)$ equals $1/p$ with $p=\dim(\mt)$.

In a nonlinear model, a design $\xi_{eG}^*$ maximizing $\phi_{eG}(\xi)$
satisfies the estimability condition \eqref{estimability} at $\mt=\mt
^0$. Indeed, $\max_{x\in\SX} [\eta(x,\mt)-\eta(x,\mt
^0) ]^2>0$
for any $\mt\neq\mt^0$ from \eqref{identifiable}, so that there
exists some $\xi\in\Xi$ such that $\phi_{eG}(\xi)>0$. Therefore,
$\phi
_{eG}(\xi_{eG}^*)>0$, and $\|\eta(\cdot,\mt)-\eta(\cdot,\mt^0)\|
_{\xi
_{eG}^*}^2=0$ implies that $\eta(x,\mt)=\eta(x,\mt^0)$ for all
$x\in\SX
$, that is, $\mt=\mt^0$ from \eqref{identifiable}. Notice that when
$\Theta$ contains an open neighborhood of $\mt^0$, then $\phi
_{eG}(\xi
)\leq1/p$ for all $\xi\in\Xi$.

Again, directional derivatives can easily be computed and
an optimal design can be obtained by linear programming when $\Theta$
and $\SX$ are both finite, or with the algorithm of Section~\ref
{S:eE-relax} when $\SX$ is finite but $\Theta$ has nonempty interior.
Note that there are now $m \times\ell$ inequality constraints in
\eqref
{CTLP}, given by
\[
\sum_{i=1}^\ell w_i
h_i\bigl(\mt^{(j)},x^{(k)}\bigr) \geq t,\qquad j=1,
\ldots,m, k=1,\ldots,\ell,
\]
where now
\[
h_{i}(\mt,x) = \frac{[\eta(x^{(i)},\mt)-\eta(x^{(i)},\mt
^0)]^2}{
[\eta(x,\mt)-\eta(x,\mt^0) ]^2}.
\]
Also note that in the algorithm of Section~\ref{S:eE-relax} we need to
construct two sequences of sets, $\Theta^{(k)}$ and $\SX^{(k)}$, with
$\Theta^{(k+1)}=\Theta^{(k)}\cup\{\mt^{(k+1)}\}$ and $\SX
^{(k+1)}=\SX
^{(k)}\cup\{\hat x^{(k+1)}\}$ at step 2, and \eqref{theta_k+1}
replaced by
\[
\bigl\{ \mt^{(k+1)}, \hat x^{(k+1)} \bigr\} =\arg\min
_{\{\mt,x\}\in\Theta
\times\SX
} \frac{\|\eta(\cdot,\mt)-\eta(\cdot,\mt^0)\|_{\xi_k}^2}{
[\eta
(x,\mt)-\eta(x,\mt^0) ]^2}
\]
with $\xi_k$ the design measure corresponding to the weights $\wb^k$.

\section{Examples}\label{S:example}

We shall use the common notation
\[
\xi= \left\{ \matrix{
 x_1 & \cdots&
x_m
\vspace*{2pt}\cr
w_1 & \cdots& w_m }
 \right\}
\]
for a discrete design measure with $m$ support points $x_i$ and such
that $\xi(\{x_i\})=w_i$, $i=1,\ldots,m$.
In the three examples considered, we indicate the values of the
parametric, intrinsic and total measure of curvatures at $\mt^0$ (for
$\sigma= 1$); see Tables~\ref{T:exPW2001}, \ref{T:c-opt-ACHJ93} and
\ref{T:KW}. They are not used for the construction of optimal designs,
and the examples illustrate the fact that they provide information on
the local behavior only (at $\mt^0$), so that a small curvature does
not mean good performance in terms of extended optimality. They are
given by
\begin{eqnarray*}
C_{\mathrm{int}}(\xi,\mt) &=& \sup_{\ub\in\mathbb{R}^p-\{\0b\}} \frac{\|
[I-P_\mt]
\sum_{i,j=1}^p u_i [\mp^2\eta(\cdot,\mt)/\mp\mt_i\,\mp\mt_j]u_j
\|_\xi} {\ub^\top \Mb(\xi,\mt)\ub},
\\
C_{\mathrm{par}}(\xi,\mt) &=& \sup_{\ub\in\mathbb{R}^p-\{\0b\}} \frac{\|
P_\mt
\sum_{i,j=1}^p u_i [\mp^2\eta(\cdot,\mt)/\mp\mt_i\,\mp\mt_j]u_j
\|_\xi} {\ub^\top \Mb(\xi,\mt)\ub},
\\
C_{\mathrm{tot}}(\xi,\mt) &=& \sup_{\ub\in\mathbb{R}^p-\{\0b\}} \frac{\|
\sum_{i,j=1}^p u_i [\mp^2\eta(\cdot,\mt)/\mp\mt_i\,\mp\mt_j]u_j
\|_\xi} {\ub^\top \Mb(\xi,\mt)\ub}
\\
& \leq& C_{\mathrm{int}}(\xi,\mt)+ C_{\mathrm{par}}(\xi,\mt),
\end{eqnarray*}
with $P_\mt$ the projector
\[
( P_\theta f ) \bigl(x'\bigr) = \frac{\mp\eta(x',\mt)}{\mp\mt
^\top } \Mb
^{-1}(\xi,\theta) \int_{\SX} \frac{\mp\eta(x,\mt)}{\mp\mt} f(x)
\xi(\D x),
\]
and correspond to the original measures of nonlinearity of \citet
{BatesW80} for $\sigma= 1$, with an adaptation to the use of a design
measure $\xi$ instead of an exact design $(x_1,\ldots,x_N)$. The
connection with the curvature arrays of \citet{BatesW80} is presented
in \citet{Pazman93}, Section~5.5; a procedure for their numerical
computation is given in \citet{BatesW80,Ratkowsky83}.

All computations are performed in Matlab on a biprocessor PC (2.5 GHz)
with 64 bits, equipped with 32 Gb RAM. Classical optimal designs ($D$-,
$E$- and $c$-optimality) are computed with
the cutting-plane method; see \citeauthor{PP2013} (\citeyear{PP2013}),
Section~9.5.3; LP problems
are solved with the simplex algorithm; we use sequential quadratic programming
for the local minimization of $H(\wb^k,\mt)$ that yields $\mt^{(k+1)}$
in \eqref{theta_k+1_V2}-(ii).

\begin{example}\label{ex2}
This example is artificial and constructed to illustrate the possible
pitfall of using a local approach (here $E$-optimal design) for
designing an experiment.
The model response is given by
\[
\eta(\xb,\mt)=\mt_1 \{\xb\}_1 + \mt_1^3
\bigl(1-\{\xb\}_1\bigr)+\mt_2 \{ \xb\}_2+\mt
_2^2\bigl(1-\{\xb\}_2\bigr), \qquad \mt=(
\mt_1,\mt_2)^\top,
\]
with $\xb\in\SX=[0,1]^2$ and $\{\xb\}_i$ denoting the $i$th component
of $\xb$. We consider local designs for $\mt^0=(1/8,1/8)^\top $.
One may notice that the set $\{\mp\eta(\xb,\mt)/\break \mp\mt
{|}_{\mt^0}\dvtx\xb\in\SX\}$ is the rectangle $[3/64,1]\times
[1/4,1]$, so that optimal
designs for any isotonic criterion function of the information matrix
$\Mb(\xi)$ are supported on the vertices $(0, 1)$, $(1, 0)$ and
$(1, 1)$ of $\SX$. The classical $D$- and $E$-optimal designs are supported
on three and two points, respectively,
\begin{eqnarray*}
\xi_{D,\mt^0}^* &\simeq&\left\{ \matrix{
\pmatrix{ 0
\cr
1} & \pmatrix{ 1
\cr
0}
& \pmatrix{1
\cr
1}
\vspace*{2pt}\cr
0.4134 & 0.3184 & 0.2682}
 \right\},\\
\xi_{E,\mt^0}^* &\simeq&\left\{\matrix{
 \pmatrix{ 0
\cr
1}
 & \pmatrix{ 1
\cr
0}
\vspace*{2pt}\cr
0.5113 & 0.4887 }
 \right\}.
\end{eqnarray*}


When only the design points $\xb_1=(0\ 1)^\top $ and $\xb_2=(1\ 0)^\top $ are
used, the parameters are only locally estimable. Indeed, the equations
in $\mt'$
\begin{eqnarray*}
\eta\bigl(\xb_1,\mt'\bigr) &=& \eta(
\xb_1,\mt),
\\
\eta\bigl(\xb_2,\mt'\bigr) &=& \eta(\xb_2,
\mt)
\end{eqnarray*}
give not only the trivial solutions $\mt'_1=\mt_1$ and $\mt'_2=\mt_2$
but also $\mt'_1$ and $\mt'_2$ as roots of two univariate polynomials
of the fifth degree (with coefficients depending on $\mt$). Since these
polynomials always admit at least one real root, at least one solution
exists for $\mt'$ that is different from $\mt$. In particular, the
vector ${\mt^0}'=(-0.9760, 1.0567)^\top $ gives approximately the same
values as $\mt^0$ for the responses at $\xb_1$ and $\xb_2$.

Direct calculations indicate that, for any $\mt$, the maximum of $\|
\eta
(\cdot,\mt)-\eta(\cdot,\mt^0)\|_\xi^2$ with respect to $\xi\in
\Xi$ is
reached for a measure supported on $(0, 0)$, $(0, 1)$, $(1, 0)$ and
$(1, 1)$. Also, the maximum of $[\eta(x,\mt)-\eta(x,\mt^0)]^2$ with
respect to $x$ is attained on the same points. We can thus restrict our
attention to the design space $\SX=\{(0, 0), (0, 1), (1, 0), (1, 1)\}$.
We take $\Theta=[-3,4]\times[-2,2]$ and use the algorithm of
Section~\ref{S:eE-relax}, with the grid $\SG^0$ of \eqref
{theta_k+1_V2}-(iii)
given by a random Latin hypercube design with $10\mbox{,}000$ points in
$[0,1]^2$ renormalized to $\Theta$ [see, e.g., \citet{tang93}], to
determine optimal designs for $\phi_{eE}(\cdot)$ and $\phi
_{eG}(\cdot)$.
When initialized with the uniform measure on the four points of $\SX$,
and with $\me=10^{-10}$, the algorithm stops after 46 and 15
iterations, respectively, requiring 0.67 s and 0.28 s in total, and
gives the designs
\begin{eqnarray*}
\xi_{eE,\mt^0}^* &\simeq& \left\{ \matrix{
\pmatrix{ 0
\cr
0} & \pmatrix{ 0
\cr
1}& \pmatrix{ 1
\cr
1}
\vspace*{2pt}\cr
0.32 & 0.197 & 0.483}
 \right\},
\\
\xi_{eG,\mt^0}^* &\simeq& \left\{ \matrix{
\pmatrix{ 0
\cr
0} & \pmatrix{ 0
\cr
1}
 & \pmatrix{ 1
\cr
0}
 & \pmatrix{ 1
\cr
1}
\vspace*{2pt}\cr
0.258 & 0.258 & 0.258 & 0.226 }
 \right\}.
\end{eqnarray*}

The performances of the designs $\xi_{D}^{*}$, $\xi_{E}^{*}$, $\xi
_{eE}^{*}$ and $\xi_{eG}^{*}$ are given in Table~\ref{T:exPW2001}. The
values $\phi_{eE}(\xi_E^*)=\phi_{eG}(\xi_E^*)=0$ indicate that
$E$-optimal design is not suitable here, the model being only locally
identifiable for $\xi_E^*$. The parametric, intrinsic and total
measures of curvature at $\mt^0$ (for $\ms^2=1$) are also indicated in
Table~\ref{T:exPW2001}.
Notice that the values of these curvature at $\mt^0$ do not reveal any
particular difficulty concerning $\xi_E^*$, but that the lack of
identifiability for this design is pointed out by the extended
optimality criteria.
\end{example}

\begin{table}
\caption{Performances of designs $\xi_{D}^{*}$, $\xi_{E}^{*}$, $\xi
_{eE}^{*}$ and $\xi_{eG}^{*}$ and curvature measures at $\mt^0$ in
Example~\protect\ref{ex2}; $\det^{1/3}=\phi_D(\xi)=\{\det[\Mb(\xi,\mt^0)]\}
^{1/3}$, $\ml
_{\min}=\phi_E(\xi)=\ml_{\min}[\Mb(\xi,\mt^0)]$. The optimal (maximum)
values of the criteria are indicated in boldface}\label{T:exPW2001}
\begin{tabular*}{\textwidth}{@{\extracolsep{\fill}}lcccccd{1.3}c@{}}
\hline
\multicolumn{1}{@{}l}{$\bolds{\xi}$} &
\multicolumn{1}{c}{$\bolds{\det^{1/3}}$} &
\multicolumn{1}{c}{$\bolds{\ml_{\min}}$} &
\multicolumn{1}{c}{$\bolds{\phi_{eE}}$} &
\multicolumn{1}{c}{$\bolds{\phi_{eG}}$} &
\multicolumn{1}{c}{$\bolds{C_{\mathrm{par}}}$} &
\multicolumn{1}{c}{$\bolds{C_{\mathrm{int}}}$} &
\multicolumn{1}{c@{}}{$\bolds{C_{\mathrm{tot}}}$}\\
\hline
$\xi_D^*$ & \textbf{0.652} & 0.273 & $3.16\cdot10^{-3}$ & 0.108 & 1.10
& 0.541 & 1.22 \\
$\xi_E^*$ & 0.625 & \textbf{0.367} & 0 & 0\phantom{000.} & 1.19 & 0 & 1.19 \\
$\xi_{eE}^*$ & 0.453 & $8.45\cdot10^{-2}$ & $\mathbf{8\bolds{.}78\bolds{\cdot}
10^{\bolds{-3}}}$ & $9.74\cdot10^{-2}$ & 3.33 & 2.69 & 4.28 \\
$\xi_{eG}^*$ & 0.540 & 0.195 & $5.68\cdot10^{-3}$ & \textbf{0.340} &
1.33 & 1.26 & 1.83 \\
\hline
\end{tabular*}
\end{table}

This example is very particular and situations where the model is
locally, but not globally, identifiable are much more common: in that
case, \eqref{identifiable} is only satisfied locally, for $\mt'$ in a
neighborhood of $\mt$, and one may refer, for example, to \citet{Walter87,WPl95}
for a precise definition and examples. The lack of global
identifiability would then not be detected by classical optimal design,
but the maximum of $\phi_{eE}(\cdot)$ and $\phi_{eG}(\cdot)$ would be
zero for $\Theta$ large enough, showing that the model is not globally
identifiable.


%

\begin{example}\label{ex3}
Consider the regression model (one-compartment with first-order
absorption input) used in \citet{AtkinsonCHJ93},
%
\begin{eqnarray}
\label{model_Ex7} \eta(x,\mt)&=&\mt_1\bigl[\exp(-\mt_2 x)-
\exp(-\mt_3 x) \bigr],
\nonumber
\\[-8pt]
\\[-8pt]
\nonumber
\mt&=&(\mt _1,\mt _2,
\mt_3)^\top, \qquad x\in\mathbb{R}^+,
\end{eqnarray}
with nominal parameters $\mt^0=(21.80, 0.05884, 4.298)^\top $. The
$D$- and $E$-optimal designs for $\mt^0$ are, respectively, given by
\begin{eqnarray*}
\xi_{D,\mt^0}^* &\simeq& \left\{\matrix{
0.229 & 1.389 & 18.42
\vspace*{2pt}\cr
1/3 & 1/3 & 1/3 }
 \right\},
\\
\xi_{E,\mt^0}^* &\simeq& \left\{\matrix{ 0.170 &
1.398 & 23.36
\vspace*{2pt}\cr
0.199 & 0.662 & 0.139}
 \right\};
\end{eqnarray*}
see \citet{AtkinsonCHJ93}.
\end{example}

We take $\Theta$ as the rectangular region $[16, 27]\times[0.03,
0.08]\times[3, 6]$ and use the algorithm of Section~\ref{S:eE-relax} to
compute an optimal design for $\phi_{eE}(\cdot)$; the grid $\SG^0$ of
\eqref{theta_k+1_V2}-(iii) is taken as a random Latin hypercube
design with $10\mbox{,}000$ points in $[0,1]^3$ renormalized to $\Theta$.
The number of iterations and computational time depend on $\ell$, the
number of elements of $\SX$. For instance, when $\SX$ is the finite set
$\{0.2, 0.4, 0.6,\ldots,24\}$ with
$\ell=120$, and the required precision $\epsilon$ equals $10^{-10}$,
the algorithm initialized at the uniform measure on the three points
0.2, 1 and 23 converges after 42 iterations in about 26 s. By refining
$\SX$ iteratively around the support points of the current optimal
design, after a few steps we obtain
\[
\xi_{eE,\mt^0}^* \simeq\left\{ \matrix{ 0.1785 &
1.520 & 20.95
\cr
0.20 & 0.66 & 0.14 }
 \right\}.
\]
A similar approach is used below for the construction of optimal
designs for $\phi_{ec}(\cdot)$ and in Example~\ref{ex4} for $\phi_{eE}(\cdot)$.
The performances of the designs $\xi_D^*$, $\xi_E^*$ and $\xi_{eE}^*$
are indicated in Table~\ref{T:c-opt-ACHJ93}. One may notice that the
design $\xi_{eE}^*$ is best or second best for $\phi_D(\cdot)$,
$\phi
_E(\cdot)$ and $\phi_{eE}(\cdot)$ among all locally optimal designs
considered.

\begin{sidewaystable}
\tabcolsep=0pt
\tablewidth=\textwidth
\caption{Performances of different designs and curvature measures at
$\mt^0$ for the model \protect\eqref{model_Ex7} with
$\mt^0=(21.80, 0.05884, 4.298)^\top $ and $\Theta=[16, 27]\times
[0.03, 0.08]\times[3, 6]$; $\det^{1/3}=\phi_D(\xi)=
\{\det[\Mb(\xi,\mt^0)]\}^{1/3}$, $\ml_{\min}=\phi_E(\xi)=\ml
_{\min}[\Mb
(\xi,\mt^0)]$. The optimal (maximum) values of the
criteria are on the main diagonal and indicated in boldface. The bottom
part of the table corresponds to the average-optimal designs of
\protect\citet{AtkinsonCHJ93}}\label{T:c-opt-ACHJ93}
\begin{tabular*}{\textwidth}{@{\extracolsep{\fill}}ld{2.3}
d{1.3}d{1.3}ccd{2.3}d{2.3}d{1.3}d{1.3}d{1.3}d{1.3}d{1.3}@{}}
\hline
\multicolumn{1}{@{}l}{$\bolds{\xi}$} &
\multicolumn{1}{c}{$\bolds{\det^{1/3}}$} &
\multicolumn{1}{c}{$\bolds{\ml_{\min}}$} &
\multicolumn{1}{c}{$\bolds{\phi_{eE}}$} &
\multicolumn{1}{c}{$\bolds{\phi_{c_1}}$}
& \multicolumn{1}{c}{$\bolds{\phi_{ec_1}}$} &
\multicolumn{1}{c}{$\bolds{\phi_{c_2}}$} &
\multicolumn{1}{c}{$\bolds{\phi_{ec_2}}$} &
\multicolumn{1}{c}{$\bolds{\phi_{c_3}}$ }
& \multicolumn{1}{c}{$\bolds{\phi_{ec_3}}$}
& \multicolumn{1}{c}{$\bolds{C_{\mathrm{par}}}$} &
\multicolumn{1}{c}{$\bolds{C_{\mathrm{int}}}$} &
\multicolumn{1}{c@{}}{$\bolds{C_{\mathrm{tot}}}$}\\
\hline
$\xi_D^*$ & \multicolumn{1}{c}{\textbf{11.74}\phantom{0}} & 0.191 & 0.178 & \multicolumn{1}{c}{$1.56\cdot10^{-4}$} &
\multicolumn{1}{c}{$6.68\cdot10^{-5}$} & 23.43 & 18.31 & 0.361 & 0.356 &
0.526 & 0 & 0.526 \\
$\xi_E^*$ & 8.82 & \multicolumn{1}{c}{\textbf{0.316}} & 0.274 & \multicolumn{1}{c}{$6.07\cdot10^{-5}$} &
\multicolumn{1}{c}{$3.08\cdot10^{-5}$} & 15.89 & 10.35 & 0.675 & 0.667 &
0.370 & 0 & 0.370 \\
$\xi_{eE}^*$ & 9.05 & 0.311 & \multicolumn{1}{c}{\textbf{0.281}} & \multicolumn{1}{c}{$6.45\cdot10^{-5}$} &
\multicolumn{1}{c}{$3.01\cdot10^{-5}$} & 16.62 & 11.03 & 0.656 & 0.644 &
0.358 & 0 & 0.358 \\
$\xi_{c_1}^*$ & 0 & 0 & 0 & \multicolumn{1}{c}{$\mathbf{4\bolds{.}56\bolds{\cdot}10^{\bolds{-4}}}$} & 0 & 0 & 0 &
0 & 0 \\
$\xi_{ec_1}^*$ & 0.757 & \multicolumn{1}{c}{$2.70\cdot10^{-3}$} & \multicolumn{1}{c}{$1.92\cdot10^{-3}$} &
\multicolumn{1}{c}{$2.26\cdot10^{-4}$} & \multicolumn{1}{c}{$\mathbf{2\bolds{.}17\bolds{\cdot}10^{\bolds{-4}}}$} & \multicolumn{1}{c}{$8.55\cdot
10^{-2}$} & \multicolumn{1}{c}{$6.12\cdot10^{-2}$} & \multicolumn{1}{c}{$1.12\cdot10^{-2}$} & \multicolumn{1}{c}{$1.09\cdot
10^{-2}$} &
6.51 & 0 & 6.51 \\
$\xi_{c_2}^*$ & 0 & 0 & 0 & 0 & 0 & \multicolumn{1}{c}{\textbf{35.55}\phantom{0}} & 0 & 0 & 0 \\
$\xi_{ec_2}^*$ & 7.86 & \multicolumn{1}{c}{$7.20\cdot10^{-2}$} & \multicolumn{1}{c}{$5.99\cdot10^{-2}$} &
\multicolumn{1}{c}{$4.55\cdot10^{-5}$} & \multicolumn{1}{c}{$1.81\cdot10^{-5}$} & 28.82 & \multicolumn{1}{c}{\textbf{27.20}\phantom{0}} &
0.157 & 0.145 &
1.12 & 0.028 & 1.12 \\
$\xi_{c_3}^*$ & 0 & 0 & 0 & 0 & 0 & 0 & 0 & \multicolumn{1}{c}{\textbf{1}\phantom{000.}} & 0 \\
$\xi_{ec_3}^*$ & 4.06 & 0.162 & 0.137 & \multicolumn{1}{c}{$9.70\cdot10^{-6}$} &
\multicolumn{1}{c}{$4.19\cdot10^{-6}$} & 6.77 & 4.36 & 0.890 & \multicolumn{1}{c}{\textbf{0.865}} &
1.11 & 0.263 & 1.14 \\[3pt]
$\xi_{AD-A}^*$ & 11.74 & 0.191 & 0.177 & \multicolumn{1}{c}{$1.56\cdot10^{-4}$} &
\multicolumn{1}{c}{$6.68\cdot10^{-5}$} & 23.53 & 18.43 & 0.360 & 0.355 & 0.522 & 0 & 0.522
\\
$\xi_{Ac_1-A}^*$ & 2.74 & \multicolumn{1}{c}{$2.07\cdot10^{-2}$} & \multicolumn{1}{c}{$1.69\cdot10^{-2}$} &
\multicolumn{1}{c}{$4.36\cdot10^{-4}$} & \multicolumn{1}{c}{$1.50\cdot10^{-4}$} & 1.12 & 0.864 & \multicolumn{1}{c}{$4.06\cdot
10^{-2}$} & \multicolumn{1}{c}{$4.01\cdot10^{-2}$} & 1.82 & 0 & 1.82 \\
$\xi_{Ac_2-A}^*$ & 6.71 & \multicolumn{1}{c}{$7.22\cdot10^{-2}$} & \multicolumn{1}{c}{$6.66\cdot10^{-2}$} &
\multicolumn{1}{c}{$2.21\cdot10^{-5}$} & \multicolumn{1}{c}{$6.67\cdot10^{-6}$} & 35.16 & 20.31 & 0.175 &
0.175 & 0.909 & 0 & 0.909 \\
$\xi_{Ac_3-A}^*$ & 3.31 & 0.118 & \multicolumn{1}{c}{$8.23\cdot10^{-2}$} & \multicolumn{1}{c}{$8.06\cdot
10^{-6}$} & \multicolumn{1}{c}{$3.86\cdot10^{-6}$} & 4.37 & 3.17 & 0.937 & 0.838 & 1.82 & 0
& 1.82 \\[3pt]
$\xi_{AD-B}^*$ & 11.08 & 0.179 & 0.159 & \multicolumn{1}{c}{$1.62\cdot10^{-4}$} &
\multicolumn{1}{c}{$6.99\cdot10^{-5}$} & 21.15 & 15.93 & 0.338 & 0.335 & 0.505 & 0.056 &
0.507 \\
$\xi_{Ac_1-B}^*$ & 2.18 & \multicolumn{1}{c}{$2.23\cdot10^{-2}$} & \multicolumn{1}{c}{$1.46\cdot10^{-2}$} &
\multicolumn{1}{c}{$2.34\cdot10^{-4}$} & \multicolumn{1}{c}{$1.56\cdot10^{-4}$} & 0.791 & 0.644 & \multicolumn{1}{c}{$5.37\cdot
10^{-2}$} & \multicolumn{1}{c}{$4.89\cdot10^{-2}$} & 2.12 & 0.133 & 2.13 \\
$\xi_{Ac_2-B}^*$ & 9.45 & 0.162 & 0.134 & \multicolumn{1}{c}{$8.07\cdot10^{-5}$} &
\multicolumn{1}{c}{$3.03\cdot10^{-5}$} & 20.05 & 16.28 & 0.385 & 0.358 & 0.753 & 0.118 &
0.761 \\
$\xi_{Ac_3-B}^*$ & 6.16 & 0.149 & \multicolumn{1}{c}{$9.92\cdot10^{-2}$} & \multicolumn{1}{c}{$5.07\cdot
10^{-5}$} & \multicolumn{1}{c}{$2.06\cdot10^{-5}$} & 6.60 & 6.13 & 0.615 & 0.587 & 1.22 &
0.256 & 1.25 \\
\hline
\end{tabular*}
\end{sidewaystable}

The intrinsic curvature is zero for $\xi_D^*$, $\xi_E^*$ and $\xi
_{eE}^*$ [since they all have $3=\dim(\mt)$ support points] and the
parametric curvatures at $\mt^0$ are rather small (the smallest one is
for $\xi_{eE}^*$). This explains that, the domain $\Theta$ being not
too large, the values of $\phi_{eE}(\xi)$ do not differ very much from
those of $\phi_E(\xi)=\ml_{\min}[\Mb(\xi,\mt^0)]$.\vadjust{\goodbreak}

Consider now the same three functions of interest as in \citet{AtkinsonCHJ93}: $ g_1(\mt)$ is the area under the curve,
\[
g_1(\mt)=\int_0^\infty\eta(x,\mt) \,\D
x = \mt_1 (1/\mt _2-1/\mt_3 );
\]
$ g_2(\mt)$ is the time to maximum concentration,
\[
g_2(\mt) =\frac{\log\mt_3 - \log\mt_2}{\mt_3-\mt_2},
\]
and $ g_3(\mt)$ is the maximum concentration,
\[
g_3(\mt) = \eta\bigl[ g_2(\mt),\mt\bigr].
\]
We shall write $\cb_i=\cb_i(\mt^0)=\mp g_i(\mt)/\mp\mt|_{\mt^0}$ and
denote $\xi_{c_i,\mt^0}$ the (locally) optimal design for $ g_i(\mt)$
which maximizes $\phi_{c_i}(\xi;\mt^0)=[\cb_i^\top \Mb^-(\xi,\mt
^0)\cb
_i]^{-1}$, for $i=1,2,3$. The $\xi_{c_i,\mt^0}^*$ are singular and are
approximately given by
\begin{eqnarray*}
\xi_{c_1,\mt^0}^* &\simeq& \left\{ \matrix{ 0.2327 &
17.63
\cr
0.0135 & 0.9865 }
 \right\},
\\
\xi_{c_2,\mt^0}^* &\simeq& \left\{ \matrix{ 0.1793 &
3.5671
\cr
0.6062 & 0.3938 }
 \right\},
\\
\xi_{c_3,\mt^0}^* &\simeq& \left\{ \matrix{ 1.0122
\cr
1 }
 \right\};
\end{eqnarray*}
see \citet{AtkinsonCHJ93}.

For each function $g_i$, we restrict the search of a design $\xi
_{ec_i}$ optimal in the sense of the criterion $\phi_{ec}(\cdot)$ to
design measures supported on the union of the supports of $\xi_{D,\mt
^0}^*$, $\xi_{E,\mt^0}^*$ and $\xi_{c_i,\mt^0}^*$.
%
%
We then obtain the following designs:
\begin{eqnarray*}
\xi_{ec_1,\mt^0}^* &\simeq& \left\{\matrix{
0.2327 & 1.389 & 23.36
\vspace*{2pt}\cr
9\cdot10^{-4} & 1.2\cdot 10^{-2} & 0.9871}
 \right\},
\\
\xi_{ec_2,\mt^0}^* &\simeq& \left\{ %
\matrix{
0.1793 & 0.229 & 3.5671 & 18.42
\cr
5.11\cdot 10^{-2} & 0.5375 & 0.3158 & 9.56\cdot10^{-2}
}
 \right\},
\\
\xi_{ec_3,\mt^0}^* &\simeq& \left\{ %
\matrix{
0.229 & 1.0122 & 1.389 & 18.42
\cr
8.42\cdot10^{-2} & 0.4867 & 0.4089 & 2.02\cdot10^{-2}
}
 \right\}.
\end{eqnarray*}

The performances of $\xi_{c_i}^*$ and $\xi_{ec_i}^*$, $i=1,\ldots,3$,
are indicated in Table~\ref{T:c-opt-ACHJ93}, together with the
curvature measures at $\mt^0$ for $\xi_{ec_i}^*$ (which are
nonsingular). For each function $ g_i$ of interest, the design $\xi
_{ec_i}^*$ performs slightly worse than $\xi_{c_i}^*$ in terms of
$c$-optimality, but contrarily to $\xi_{c_i}^*$, it allows us to
estimate the three parameters $\mt$ and guarantees good estimability
properties for $ g_i(\mt)$ for all $\mt\in\Theta$. Notice that, apart
from the $c$-optimality criteria $\phi_{c_i}(\cdot)$, all criteria
considered take the value 0 at the $c$-optimal designs $\xi_{c_i}^*$.
The construction of an optimal design for $\phi_{ec}(\cdot)$ thus forms
an efficient method to circumvent the difficulties caused by singular
$c$-optimal design in nonlinear models;
see \citet{PP2013}, Chapters~3 and 5.
One may also refer to \citet{PMSlovaca-09} for alternative
approaches for the regularization of singular $c$-optimal designs.

We conclude this example with a comparison with the average-optimal
designs of \citet{AtkinsonCHJ93} that aim at taking uncertainty on
$\mt
^0$ into account. Consider a prior distribution $\pi(\cdot)$ on the two
components of $\mt$ that intervene nonlinearly in $\eta(x,\mt)$, and
let $\Ex_\pi\{\cdot\}$ denote the expectation for $\pi(\cdot)$.
\citet{AtkinsonCHJ93} indicate that when $\pi$ equals $\pi_A$ uniform
on $[\mt_2^0-0.01,\mt_2^0+0.01]\times[\mt_3^0-1, \mt_3^0+1]$, the
design that maximizes $\Ex_\pi\{\log\det[\Mb(\xi,\mt)]\}$ is
\[
\xi_{AD-A}^* \simeq\left\{ \matrix{ 0.2288 &
1.4170 & 18.4513
\cr
1/3 & 1/3 & 1/3 }
 \right\},
\]
and the designs that minimize $\Ex_\pi\{ \cb_i^\top (\mt)\Mb
^-(\xi,\mt)\cb
_i(\mt) \}$, $i=1,2,3$, are
\begin{eqnarray*}
\xi_{Ac_1-A}^* &\simeq& \left\{\matrix{ 0.2449 &
1.4950 & 18.4903
\cr
0.0129 & 0.0387 & 0.9484 }
 \right\},
\\
\xi_{Ac_2-A}^* &\simeq& \left\{\matrix{ 0.1829 &
2.4639 & 8.8542
\cr
0.6023 & 0.2979 & 0.0998 }
 \right\},
\\
\xi_{Ac_3-A}^* &\simeq& \left\{ %
\matrix{ 0.3608 &
1.1446 & 20.9218
\cr
0.0730 & 0.9094 & 0.0176 }
 \right\}.
\end{eqnarray*}
When $\pi$ equals $\pi_B$ uniform on $[\mt_2^0-0.04,\mt
_2^0+0.04]\times
[\mt_3^0-4, \mt_3^0+4]$, the average-optimal designs are
\begin{eqnarray*}
\xi_{AD-B}^* &\simeq& \left\{\matrix{ 0.2034
& 1.1967 & 2.8323 & 7.8229 & 20.1899
\cr
0.2870 & 0.2327 & 0.1004 & 0.0678 & 0.3120 }
 \right\},
\\
\xi_{Ac_1-B}^* &\simeq& \left\{ %
\matrix{ 0.2909 &
1.7269 & 13.0961 & 39.58
\cr
0.0089 & 0.0365 & 0.2570 & 0.6976 }
 \right\},
\\
\xi_{Ac_2-B}^* &\simeq& \left\{\matrix{ 0.2513
& 0.9383 & 2.7558 & 8.8381 & 26.6564
\cr
0.2914 & 0.2854 & 0.1468 & 0.2174 & 0.0590 }
 \right\},
\\
\xi_{Ac_3-B}^* &\simeq& \left\{ %
\matrix{ 0.3696
& 1.1383 & 2.4370 & 6.0691 & 24.0831
\cr
0.0971 & 0.3584 & 0.3169 & 0.1634 & 0.0641}
 \right\}.
\end{eqnarray*}
Their performances are indicated in the bottom part of Table~\ref
{T:c-opt-ACHJ93}.
Notice that the average-optimal designs for the vague prior $\pi_B$ are
supported on more than three points and thus allow model checking; this
is the case too for the two designs $\xi_{ec_2}^*$ and $\xi_{ec_3}^*$.
However, contrary to average-optimal design, the number of support
points of optimal designs for extended optimality criteria does
not\vadjust{\goodbreak}
seem to increase with uncertainty measured by the size of $\Theta$: for
instance, when $\Theta=[\mt_1^0-5,\mt_1^0+5]\times[\mt
_2^0-0.04,\mt
_2^0+0.04]\times[\mt_3^0-4, \mt_3^0+4]$, the optimal design for
$\phi
_{eE}(\cdot)$ is still supported on three points, approximately 0.1565,
1.552 and 19.73, receiving weights 0.268, 0.588 and 0.144, respectively.

\begin{figure}

\includegraphics{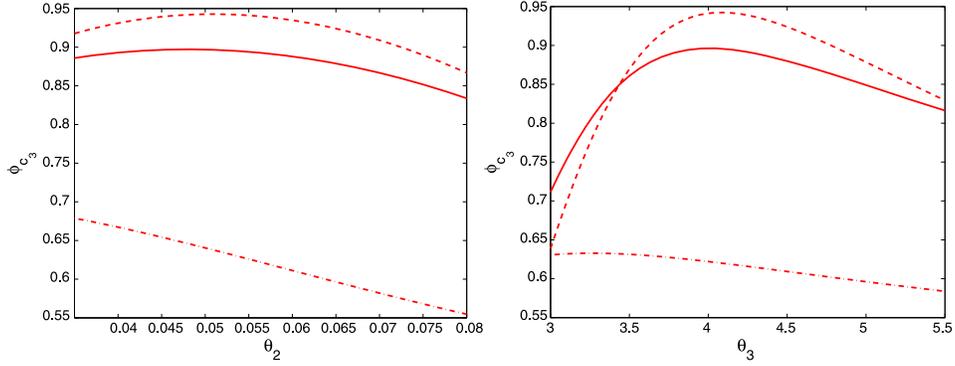}

\caption{$\phi_{c_3}(\xi;\mt)$ as a function of $\mt$, for $\mt
_1=\mt
_1^0$ and $\mt_3=\mt_3^0$ (\emph{left})
and $\mt_1=\mt_1^0$, $\mt_2=\mt_2^0$ (\emph{right}); $\xi=\xi
_{Ac_3-A}^*$ in \emph{dashed line}, $\xi=\xi_{Ac_3-B}^*$ in
\emph{dash--dotted line} and $\xi=\xi_{ec_3,\mt^0}^*$ in \emph
{solid line}.}
\label{F:Ex3}
\end{figure}

All average-optimal designs considered yield reasonably small
curvatures at $\mt^0$, although larger than those for $\xi_{E,\mt^0}^*$
and $\xi_{eE,\mt^0}^*$. The performances of $\xi_{AD-A}^*$ and $\xi
_{AD-B}^*$ are close to those of $\xi_{D,\mt^0}^*$, and the most
interesting features concern designs for estimation of functions of
interest $g_i(\mt)$. The designs $\xi_{c_i,\mt^0}^*$ cannot be used if
$\mt\neq\mt^0$ and are thus useless in practice. The average-optimal
designs $\xi_{Ac_i-B}^*$ perform significantly worse than $\xi
_{ec_i,\mt
^0}^*$ in terms of $\phi_{ec_i}(\cdot)$ for $i=1,2$ and~3 and in terms
of $\phi_{c_i}(\cdot)$ for $i=2$ and 3. On the other hand, the designs
$\xi_{Ac_i-A}^*$, constructed for the precise prior $\pi_A$, perform
significantly better than $\xi_{ec_i,\mt^0}^*$ in terms of $\phi
_{c_i}(\cdot)$ for all $i$. Figure~\ref{F:Ex3} presents $\phi
_{c_3}(\xi;\mt)$ as a function of $\mt$, for the three designs $\xi_{Ac_3-A}^*$
(\emph{dashed line}), $\xi_{Ac_3-B}^*$ (\emph{dash--dotted line}) and
$\xi_{ec_3,\mt^0}^*$ (\emph{solid line}), when $\mt_1=\mt_1^0$,
$\mt
_3=\mt_3^0$ (\emph{left}) and $\mt_1=\mt_1^0$, $\mt_2=\mt_2^0$
(\emph{right}).
Note that the projection on the last two components of $\mt$ of the set
$\Theta$ used for extended $c$-optimality is intermediate between the
supports of $\pi_A$ and $\pi_B$. Although average-optimal designs
$\xi
_{Ac_i-A,B}^*$ and extended-optimal designs $\xi_{ec_i,\mt^0}^*$ pursue
different objectives, the example indicates that they show some
resemblance in terms of precision of estimation of $g_i(\mt)$. The
situation would be totally different in absence of global
identifiability for $g_i(\mt)$, a problem that would not be detected by
average-optimal designs; see the discussion at the end of Example~\ref{ex2}.

\begin{example}\label{ex4}
For the same regression model \eqref{model_Ex7}, we change the value of
$\mt^0$ and the set $\Theta$ and take $\mt^0=(0.773, 0.214,
2.09)^\top $ and $\Theta=[0,5]\times[0,5]\times[0,5]$, the
values\vadjust{\goodbreak}
used by
\citet{KiefferW98}. With these values, from an investigation based on
interval analysis, the authors report that for the 16-point design
\[
\xi_0= \left\{\matrix{ %
1 & 2 & \cdots& 16
\vspace*{2pt}\cr
1/16 & 1/16 & \cdots& 1/16}\right\}
\]
and the observations $\yb$ given in their Table~13.1, the LS criterion
$\|\yb-\eta_X(\mt)\|^2$ has a global minimizer (the value we have taken
here for $\mt^0$) and two other local minimizers in $\Theta$. The $D$-
and $E$-optimal designs for $\mt^0$ are now given by
\begin{eqnarray*}
\xi_{D,\mt^0}^* &\simeq& \left\{\matrix{ %
0.42 &
1.82 & 6.80
\vspace*{2pt}\cr
1/3 & 1/3 & 1/3 }
 \right\},
\\
\xi_{E,\mt^0}^* &\simeq& \left\{ \matrix{
0.29 &
1.83 & 9.0
\vspace*{2pt}\cr
0.4424 & 0.3318 & 0.2258 }
 \right\}.
\end{eqnarray*}
Using the same approach as in Example~\ref{ex3}, with the grid $\SG^0$ of
\eqref
{theta_k+1_V2}-(iii) obtained from a random Latin hypercube design
with $10\mbox{,}000$ points in $\Theta$, we obtain
\[
\xi_{eE,\mt^0}^* \simeq\left\{ \matrix{
0.38 & 2.26
& 7.91
\vspace*{2pt}\cr
0.314 & 0.226 & 0.460}
 \right\}.
\]

To compute an optimal design for $\phi_{eG}(\cdot)$, we consider the
design space $\SX=\{0,0.1,0.2,\ldots,16\}$ (with 161 points) and use
the algorithm of Section~\ref{S:eE-relax} with the grid $\SG^0$ of
\eqref{theta_k+1_V2}-(iii) taken as a random Latin hypercube design
with $10^5$ points. The same design space is used to evaluate $\phi
_{eG}(\cdot)$ for the four designs above. For $\epsilon= 10^{-10}$,
the algorithm initialized at the uniform measure on $\SX$ converges
after 34 iterations in about 52 s and gives
\[
\xi_{eG,\mt^0}^* \simeq\left\{ \matrix{
 0.4 & 1.9 &
5.3 & 16
\cr
0.278 & 0.258 & 0.244 & 0.22 }
 \right\}.
\]
\end{example}

\begin{table}
\caption{Performances of different designs and curvature measures at
$\mt^0$ for the model \protect\eqref{model_Ex7} with
$\mt^0=(0.773, 0.214, 2.09)^\top $ and $\Theta=[0,5]^3$; $\det
^{1/3}=\phi_D(\xi)=\{\det[\Mb(\xi,\mt^0)]\}^{1/3}$, $\ml_{\min
}=\phi
_E(\xi)=\ml_{\min}[\Mb(\xi,\mt^0)]$.
The optimal (maximum) values of the criteria are indicated in boldface}%
\label{T:KW}
\begin{tabular*}{\textwidth}{@{\extracolsep{\fill}}lccccd{3.1}d{2.2}d{3.1}@{}}
\hline
$\bolds{\xi}$ & $\bolds{\det^{1/3}}$ & $\bolds{\ml_{\min}}$ &
$\bolds{\phi_{eE}}$ & $\bolds{\phi_{eG}}$ &
\multicolumn{1}{c}{$\bolds{C_{\mathrm{par}}}$} &
\multicolumn{1}{c}{$\bolds{C_{\mathrm{int}}}$} &
\multicolumn{1}{c@{}}{$\bolds{C_{\mathrm{tot}}}$}\\
\hline
$\xi_0$ & $1.85\cdot10^{-2}$ & $1.92\cdot10^{-4}$ & $2.28\cdot
10^{-5}$ & $5.66\cdot10^{-3}$ & 180.7 & 15.73 & 181.3\\
$\xi_D^*$ & $\mathbf{5\bolds{.}19\bolds{\cdot}10^{\bolds{-2}}}$ & $1.69\cdot10^{-3}$ &
$2.64\cdot10^{-4}$ & $6.70\cdot10^{-2}$ & 58.0 & 0 & 58.0 \\
$\xi_E^*$ & $4.51\cdot10^{-2}$ & $\mathbf{2\bolds{.}04\bolds{\cdot}10^{\bolds{-3}}}$ &
$1.32\cdot10^{-4}$ & $7.95\cdot10^{-2}$ & 50.7 & 0 & 50.7 \\
$\xi_{eE}^*$ & $4.73\cdot10^{-2}$ & $1.53\cdot10^{-3}$ & $\mathbf
{2\bolds{.}92\bolds{\cdot}10^{\bolds{-4}}}$ & 0.114 & 54.6 & 0 & 54.6 \\
$\xi_{eG}^*$ & $4.11\cdot10^{-2}$ & $1.31\cdot10^{-3}$ & $1.69\cdot
10^{-4}$ & \textbf{0.244} & 69.7 & 10.7 & 69.9 \\
\hline
\end{tabular*}
\end{table}

The performances and curvature measures at $\mt^0$ of $\xi_0$, $\xi
_D^*$, $\xi_E^*$, $\xi_{eE}^*$ and $\xi_{eG}^*$ are given in
Table~\ref{T:KW}. The large intrinsic curvature for $\xi_0$,
associated with the
small values of $\phi_{eE}(\xi^0)$ and $\phi_{eG}(\xi^0)$, explains the
presence of local minimizers for the LS criterion, and thus the
possible difficulties for the estimation of $\mt$. The values of $\phi
_{eE}(\cdot)$ and $\phi_{eG}(\cdot)$ reported in the table indicate
that $\xi_D^*$, $\xi_E^*$, $\xi_{eE}^*$ or $\xi_{eG}^*$ would have
caused less difficulties.

\section{Further extensions and developments}\label{S:Extensions}

\subsection{An extra tuning parameter for a smooth transition to usual
design criteria}\label{S:extra K}


The criterion $\phi_{eE}(\xi;\mt^0)$ can be written as
%
\begin{eqnarray}
\label{eE-equiv-K=0} &&\phi_{eE}\bigl(\xi;\mt^0\bigr)= \max\bigl\{
\ma\in\mathbb{R}\dvtx\bigl\|\eta(\cdot,\mt )-\eta \bigl(\cdot,\mt^0\bigr)
\bigr\|_\xi^2 \geq\ma\bigl\|\mt-\mt^0\bigr\|^2,
\nonumber
\\[-8pt]
\\[-8pt]
\nonumber
&&\hspace*{218pt}\mbox{for all } \mt\in \Theta\bigr\}.
\end{eqnarray}
Instead of giving the same importance to all $\mt$ whatever their
distance to $\mt^0$, one may wish to introduce a saturation and reduce
the importance given to those $\mt$ very far from $\mt^0$, that is, consider
%
\begin{eqnarray}
&&\phi_{eE|K}\bigl(\xi;\mt^0\bigr) = \max \biggl\{\ma\in
\mathbb{R}\dvtx\bigl\| \eta (\cdot,\mt )-\eta\bigl(\cdot,\mt^0\bigr)
\bigr\|_\xi^2 \geq\ma\frac{\|\mt-\mt^0\|^2}{1+K
\|\mt
-\mt^0\|^2},
\nonumber
\\[-8pt]
\\[-8pt]
\nonumber
&&  \hspace*{255pt}\mbox{ for all } \mt\in\Theta \biggr\} \label
{eE-equiv-K}
\end{eqnarray}
for some $K\geq0$. Equivalently, $\phi_{eE|K}(\xi;\mt^0)=\min_{\mt
\in
\Theta} H_{E|K}(\xi,\mt)$, with
\[
H_{E|K}(\xi,\mt)= \bigl\|\eta(\cdot,\mt)-\eta\bigl(\cdot,\mt^0
\bigr)\bigr\|_\xi^2 \biggl[K+\frac{1}{\|\mt-\mt^0\|^{2}} \biggr].
\]
%
%
As in Section~\ref{S:eE-def}, we obtain $\phi_{eE|K}(\xi)=\ml_{\min
}[\Mb
(\xi)]$ in a linear model and, for a nonlinear model with $\Theta=\SB
(\mt^0,\rho)$, $\lim_{\rho\ra0} \phi_{eE|K}(\xi;\mt^0) = \ml
_{\min}[\Mb
(\xi,\mt^0)]$ for any $K\geq0$. Moreover, in a nonlinear model with no
overlapping $\phi_{eE|K}(\xi;\mt^0)$ can be made arbitrarily close to
$\ml_{\min}[\Mb(\xi,\mt^0)]$ by choosing $K$ large enough, whereas
choosing $K$ not too large ensures some protection against
$\|\eta(\cdot,\theta)-\eta(\cdot,\theta^0)\|_\xi$
being small for some $\mt$ far from $\mt^0$. Also,
properties of $\phi_{eE}(\cdot)$ such as concavity, positive
homogeneity, existence of directional derivatives; see Section~\ref
{S:eE-prop}, remain valid for $\phi_{eE|K}(\cdot)$, for any $K\geq0$.
The maximization of $\phi_{eE|K}(\cdot)$ forms a LP problem when both
$\SX$ and $\Theta$ are finite (see Section~\ref{S:eE-LP}) and a
relaxation procedure (cutting-plane method) can be used when $\Theta$
is a compact subset of $\mathbb{R}^p$; see Section~\ref{S:eE-relax}.

A similar approach can be used with extended $c$- and $G$-optimality,
which gives $\phi_{ec|K}(\xi)= \min_{\mt'\in\Theta} H_{c|K}(\xi,\mt')$ with
\[
H_{c|K}(\xi,\mt)= \bigl\|\eta(\cdot,\mt)-\eta\bigl(\cdot,\mt^0
\bigr)\bigr\|_\xi^2 \biggl[K+\frac{1}{| g(\mt)- g(\mt^0)|^{2}} \biggr]
\]
and
\[
\phi_{eG|K}(\xi) = \min_{\mt\in\Theta} \biggl\{ \bigl\|\eta(\cdot,
\mt)-\eta \bigl(\cdot,\mt^0\bigr)\bigr\|_\xi^2
\biggl[K+\frac{1}{\max_{x\in\SX}
[\eta(x,\mt
)-\eta(x,\mt^0) ]^2} \biggr] \biggr\},
\]
for $K$ a positive constant.

\subsection{Worst-case extended optimality criteria}\label
{S:Maximin-extended-Eopt}

The criterion defined by
\[
\phi_{MeE}(\xi)=\min_{\mt^0\in\Theta}\phi_{eE}\bigl(
\xi;\mt^0\bigr)=\min_{(\mt,\mt
^0)\in\Theta\times\Theta} H_E\bigl(
\xi,\mt;\mt^0\bigr),
\]
see \eqref{phi_eE}, \eqref{H_E}, accounts for the global behavior of
$\eta(\cdot,\mt)$ for $\mt\in\Theta$
and obliterates the dependence on $\mt^0$ that is present in $\phi
_{eE}(\xi;\mt^0)$.
The situation is similar to that in Section~\ref{S:extended-Eopt},
excepted that we consider now the minimum of $H_E$ with respect to two
vectors $\mt$ and $\mt^0$ in $\Theta\times\Theta$. All the developments
in Section~\ref{S:extended-Eopt} obviously remain valid (concavity,
existence of directional derivative, 
etc.), including the algorithmic solutions of Sections~\ref{S:eE-LP}
and \ref{S:eE-relax}.
The same is true for the worst-case versions of $\phi_{ec}(\cdot)$ and
$\phi_{eG}(\cdot)$, respectively, defined by
$\phi_{Mec}(\xi)=\min_{(\mt,\mt^0)\in\Theta\times\Theta}
H_c(\xi,\mt;\mt
^0)$, see \eqref{H_c}, and by
$\phi_{MeG}(\xi) = \min_{(\mt,\mt^0)\in\Theta\times\Theta} \{
\|\eta
(\cdot,\mt)-\eta(\cdot,\mt^0)\|_\xi^2/ \max_{x\in\SX}
[\eta(x,\mt
)-\eta(x,\mt^0) ]^2 \} $,
and for the worst-case versions of the extensions of previous section
that include an additional tuning parameter $K$.

Note that the criterion $\phi_{MeE}(\cdot)$ may direct attention to a
particularly pessimistic situation. Indeed, for $\Theta$ a compact set
with nonempty interior and $\mu$ the Lebesgue measure on $\Theta$, one
may have $\min_{\mt^0\in\Theta}\phi_{eE}(\xi;\mt^0)=0$ for all designs
$\xi$ although $\mu\{\mt^0\in\Theta\dvtx \phi_{eE}(\xi';\mt
^0)>0\}
=1$ for
some design $\xi'$. This corresponds to a situation where the model is
structurally identifiable, in the sense that the property \eqref
{identifiable} is generic but is possibly false for $\mt$ in a subset
of zero measure; see, for example, \citet{Walter87}. Example~\ref{ex2} gives an
illustration.

\setcounter{example}{1}
\begin{example}[(Continued)]
When the three polynomial equations $\mt_1'-{\mt_1'}^3=\mt_1-\mt_1^3$,
$\mt_2'-{\mt_2'}^2=\mt_2-\mt_2^2$, ${\mt_1'}^3+{\mt_2'}^2=\mt
_1^3+\mt
_2^2$ are satisfied, then $\eta(\xb,\mt')=\eta(\xb,\mt)$ for all
$\xb$.
Since these equations have solutions $\mt'\neq\mt$ in $\Theta\times
\Theta$, $\phi_{MeE}(\xi)=0$ for all $\xi\in\Xi$. On the other hand,
$\max_{\xi\in\Xi}\phi_{eE}(\xi;\mt^0)>0$ w.p.1. when $\mt^0$ is
randomly drawn with a probability measure having a density with respect
to the Lebesgue measure on $\Theta$.

In a less pessimistic version of worst-case extended $E$-optimality, we
may thus consider a finite set $\Theta^0 \subset\Theta$ for $\mt^0$,
obtained for instance by random sampling in~$\Theta$, and maximize
$\min_{\mt^0\in\Theta^0}\phi_{eE}(\xi;\mt^0)$.
\end{example}


%

\section{Conclusions}\label{S:conclusions}

Two essential ideas have been presented. First, classical optimality
criteria can be extended in a mathematically consistent way to criteria
that preserve a nonlinear model against overlapping, and at the same
time retain the main features of classical criteria, especially
concavity. Moreover, they coincide with their classical counterpart for
linear models. Second, designs that are nearly optimal for those
extended criteria can be obtained by standard linear programming
solvers, supposing that the approximation of the feasible parameter
space $\Theta$ by a finite set is acceptable. A relaxation method,
equivalent to the cutting-plane algorithm, can be used when $\Theta$ is
a compact set with nonempty interior. Linear constraints on the design
can easily be taken into account. As a by-product, this also provides
simple algorithmic procedures for the determination of $E$-, $c$- or
$G$-optimal designs in linear models with linear cost constraints.

As it is usually the case for optimal design in nonlinear models, the
extended-optimality criteria are local and depend on a guessed value
$\mt^0$ for the model parameters. However, the construction of a
globalized, worst-case version enjoying the same properties is
straightforward (Section~\ref{S:Maximin-extended-Eopt}).

Finally, we recommend the following general procedure for optimal
design in nonlinear regression. (i) Choose a parameter space $\Theta$
corresponding to the domain of interest for $\mt$, select (e.g.,
randomly) a finite subset $\Theta^0$ in the interior of $\Theta$;
(ii) for each $\mt^0$ in $\Theta^0$ compute an optimal design $\xi
_{eE,\mt^0}^*$ maximizing $\phi_{eE}(\xi;\mt^0)$ and a $E$-optimal
design $\xi_{E,\mt^0}^*$ maximizing $\phi_E(\xi;\mt^0)=\ml_{\min
}\Mb(\xi,\mt^0)$; (iii) if $\phi_{eE}(\xi_{eE,\mt^0}^*;\mt^0)$ is close
enough to $\phi_{E}(\xi_{E,\mt^0}^*;\mt^0)$ for all $\mt^0$ in
$\Theta
^0$, one may consider that the risk of overlapping, or lack of
identifiability in $\Theta$, is weak and classical optimal design that
focuses on the precision of estimation can be used; otherwise, a design
that maximizes $\min_{\mt^0\in\Theta^0} \phi_{eE}(\xi;\mt^0)$
should be
preferred. When the extended $G$-optimality criterion $\phi_{eG}(\cdot;\mt^0)$ is substituted for $\phi_{eE}(\cdot;\mt^0)$, the
comparison in
(iii) should be between $\phi_{eG}(\xi_{eG,\mt^0}^*;\mt^0)$ and
$1/\dim(\mt)$, see Section~\ref{S:E-G}. Extended $c$-optimality can be
used when one is interested in estimating a (nonlinear) function of
$\mt
$, the comparison in (iii) should then be with $c$-optimality.

\section*{Acknowledgements}
The authors thank the referees for useful comments that helped to
significantly improve the paper.




%

%



\printaddresses
\end{document}